\documentclass[smallextended]{svjour3}

\usepackage{ifthen}

\newboolean{AppendixOutsourcing}
\setboolean{AppendixOutsourcing}{true}

\newboolean{ResultTablesBlack}
\setboolean{ResultTablesBlack}{true}

\newboolean{ResultTablesNewOrder}
\setboolean{ResultTablesNewOrder}{true}

\newboolean{ResultTablesHigherDimensionalStrongTermination}
\setboolean{ResultTablesHigherDimensionalStrongTermination}{true}

\newboolean{NurQCQPversion}
\setboolean{NurQCQPversion}{true}

\usepackage[ngerman,english]{babel}
\usepackage[T1]{fontenc}
\usepackage[latin1]{inputenc}
\usepackage[numbers,comma,sort&compress]{natbib}
\usepackage{verbatim} 
\usepackage{multirow}
\usepackage{hyperref}

\usepackage{pdflscape} 
\usepackage{longtable}

\newif\ifpdf
\ifx\pdfoutput\undefined
\usepackage{graphicx}
\pdffalse
\else
\pdftrue
\usepackage{graphicx}
\usepackage{epstopdf}
\fi

\usepackage{fancyhdr}

\usepackage{amsmath,amsfonts,amssymb,amscd,xspace}
\newtheorem{algorithm}[theorem]{Algorithm}     
\newtheorem{notation}[theorem]{Notation}     
\usepackage[centerlast,small,sc]{caption}

\usepackage[numbers,comma,sort&compress]{natbib}
\usepackage{url}
\newcommand{\komma}{\textnormal{~,}}
\newcommand{\punkt}{\textnormal{~.}}

\newcommand{\xor}{\textnormal{or}}

\newcommand{\rp}{\textnormal{rp}}
\newcommand{\NforMPBNGC}{\textnormal{Nb}}
\newcommand{\NforSolvOpt}{\textnormal{Nc}}
\newcommand{\TotalTime}{\textnormal{$\textnormal{t}_{1}$}}
\newcommand{\PartTime}{\textnormal{$\textnormal{t}_{2}$}}

\newcommand{\Fscal}{c}

\renewcommand{\forall}{\textnormal{for all }}

\newcommand{\zeroVector}[1]{0}

\newcommand{\zeroMatrix}[1]{0}

\newcommand{\zeroVectorI}[1]{0}

\newcommand{\zeroVectorIT}[1]{0}

\newcommand{\Sym}[1]{\mathbb{R}_{\mathrm{sym}}^{#1\times#1}}

\newcommand{\Triu}[1]{\mathbb{R}_{\mathrm{triu}}^{#1\times#1}}

\newcommand{\Rpos}{\mathbb{R}_{\geq0}}

\newcommand{\subhessian}[2]{\partial^2{#1}{(#2)}}

\usepackage{color} 
\definecolor{gold}{rgb}{0.85,.66,0}
\definecolor{turquoise}{rgb}{0.251,0.878,0.816}

\definecolor{orange}{rgb}{1,0.5,0}   
\definecolor{violet}{rgb}{0.5,0,0.5} 

\ifthenelse{\boolean{ResultTablesBlack}}
{

}
{

}

\newcommand{\GenaueAngabeANone}{}
\newcommand{\GenaueAngabeHD}{}
\newcommand{\GenaueAngabeMIFFLIN}{}
\newcommand{\GenaueAngabeCertificate}{} 

\newcommand{\GenaueAngabeConvergenceTheory}{p.~7,~3.1~Theoretical basics}   

\newcommand{\GenaueAngabeNEWone}{p.~116, Subsection~4.3.2}  
\newcommand{\GenaueAngabeNEWtwo}{p.~123, Equation~(4.57)}  

\newcommand{\GenaueAngabeOne}{}
\newcommand{\GenaueAngabeTwo}{}
\newcommand{\GenaueAngabeThree}{}
\newcommand{\GenaueAngabeFour}{p.~9,~Remark~3.6}                            
\newcommand{\GenaueAngabeFive}{p.~10,~Example~3.7}                          
\newcommand{\GenaueAngabeSix}{p.~9,~Optimization~problems~(55)~and~(56)}    
\newcommand{\GenaueAngabeSeven}{p.~5,~Equation~(35)}                        
\newcommand{\GenaueAngabeEight}{p.~5,~Equation~(35)}                        
\newcommand{\GenaueAngabeTen}{p.~9,~Optimization~problems~(55)~and~(56)}    
\newcommand{\GenaueAngabeEleven}{p.~19,~Remark~3.16}                        
\newcommand{\GenaueAngabeFifteen}{p.~5,~Equation~(35)}                      
\newcommand{\GenaueAngabeSixteen}{p.~9,~Optimization~problem~(55)~and~(56)} 
\newcommand{\GenaueAngabeSeventeen}{p.~5,~Equation~(35)}                    

\newcommand{\refh}[1]{\textnormal{(\ref{#1})}}   

\def\fullTitle{{A feasible second order bundle algorithm for nonsmooth nonconvex optimization problems with inequality constraints: II. Implementation and numerical results}}

\def\AuthorOne{{Hannes Fendl}}
\def\AuthorThree{{Hermann Schichl}}

\ifpdf
\hypersetup{
    pdftitle  = {\fullTitle},
    pdfauthor = {\AuthorOne and \AuthorThree},
}
\fi

\ifpdf
	\usepackage[color]{changebar}
	\newcommand{\cbstartDVI}{\cbstart}
	\newcommand{\cbendDVI}{\cbend}
	
	\nochangebars  
	
  \cbcolor{blue}
	\newcommand{\cbstartMIFFLIN}{\cbstart}  
	\newcommand{\cbendMIFFLIN}{\cbend}      
	
	\renewcommand{\cbstartDVI}{{}}          
	\renewcommand{\cbendDVI}{{}}            
\else
  \newcommand{\cbstartDVI}{{}}
	\newcommand{\cbendDVI}{{}}
	
	\newcommand{\cbstartMIFFLIN}{{}}  
	\newcommand{\cbendMIFFLIN}{{}}    
\fi

\newlength{\widebarargwidth}
\newlength{\widebarwidth}
\newlength{\widebarargheight}
\newlength{\widebarargdepth}
\DeclareRobustCommand{\widebar}[1]{%
  \settowidth{\widebarargwidth}{\ensuremath{#1}}%
  \settoheight{\widebarargheight}{\ensuremath{#1}}%
  \settodepth{\widebarargdepth}{\ensuremath{#1}}%
  \addtolength{\widebarargwidth}{-0.2\widebarargheight}%
  \addtolength{\widebarargwidth}{-0.2\widebarargdepth}%
  \makebox[0pt][l]{\addtolength{\widebarargheight}{0.3ex}%
    \hspace{0.2\widebarargheight}%
    \hspace{0.2\widebarargdepth}%
    \hspace{0.5\widebarargwidth}%
    \setlength{\widebarwidth}{0.6\widebarargwidth}%
    \addtolength{\widebarwidth}{0.3ex}%
    \makebox[0pt][c]{\rule[\widebarargheight]{\widebarwidth}{0.1ex}}}%
  {#1}}

\newcommand{\emptyh}[1]{}

\begin{document}

\title{\fullTitle}
\author{\AuthorOne\thanks{This research was supported by the Austrian Science Found (FWF) Grant Nr.~P22239-N13.} \and \AuthorThree}
\institute{Faculty of Mathematics, University of Vienna, Austria\\
  Oskar-Morgenstern-Pl.~1, A-1090 Wien, Austria\\
\email{hermann.schichl@univie.ac.at}}


\maketitle

\begin{abstract}
This paper presents a concrete implementation of the
feasible second order bundle algorithm for nonsmooth, nonconvex
optimization problems with inequality constraints \cite{HannesPaperB}. It computes
the search direction by solving a convex quadratically constrained
quadratic program. Furthermore, certain versions of the search
direction problem are discussed and the applicability of this approach
is justified numerically by using different solvers for the
computation of the search direction. Finally, the good performance of
the second order bundle algorithm is demonstrated by comparison with
test results of other solvers on examples of the Hock-Schittkowski
collection, on custom examples that arise in the context of finding
exclusion boxes for quadratic constraint satisfaction problems, and on
higher dimensional piecewise quadratic examples.
\end{abstract}
\keywords{Nonsmooth optimization, nonconvex optimization, bundle method}
\subclass{90C56, 49M37, 90C30}


\section{Introduction}
Nonsmooth optimization addresses to solve the optimization problem
\begin{equation}
\begin{split}
&\min{f(x)}\\
&\textnormal{ s.t. }F_i(x)\leq0~~~\forall i=1,\dots,m\komma
\label{AOB:OptProb}
\end{split}
\end{equation}
where $f,F_i:\mathbb{R}^n\longrightarrow\mathbb{R}$ are locally Lipschitz continuous.
Since $F_i(x)\leq0$ for all $i=1,\dots,m$ if and only if
$F(x):=\max_{i=1,\dots,m}{\Fscal_iF_i(x)}\leq0$ with constants
$\Fscal_i>0$ and since $F$ is still locally Lipschitz continuous (cf.,
e.g., \citet[p.~969,~Theorem~6~(a)]{Mifflin}), we can always assume
$m=1$ in (\ref{AOB:OptProb}).
Therefore w.l.o.g.\ we always consider the nonsmooth optimization
problem with a single nonsmooth constraint
\begin{equation}
\begin{split}
&\min{f(x)}\\
&\textnormal{ s.t. }F(x)\leq0\komma
\label{BundleSQP:OptProblem}
\end{split}
\end{equation}
where $F:\mathbb{R}^n\longrightarrow\mathbb{R}$ is locally Lipschitz
continuous.

Since locally Lipschitz continuous functions are differentiable almost everywhere, both $f$ and $F$ may have kinks and therefore
already
the attempt to solve an unconstrained nonsmooth optimization problem by a smooth solver (e.g., by a line search algorithm or by a trust region method) by just replacing the gradient by a subgradient, fails in general (cf., e.g., \citet[p.~461-462]{Zowe}):
If
$g$ is an element of the subdifferential $\partial f(x)$,
then the search direction $-g$ does not need to be a direction of descent (contrary to the behavior of the gradient of a differentiable function).
Furthermore,
it
can happen that $\lbrace x_k\rbrace$ converges towards a minimizer $\hat{x}$, although the sequence of gradients $\lbrace\nabla f(x_k)\rbrace$ does not converge towards $\zeroVector{n}$ and therefore we cannot identify $\hat{x}$ as a minimizer.
Moreover,
it
can happen that $\lbrace x_k\rbrace$ converges towards a point $\hat{x}$, but $\hat{x}$ is not stationary for $f$.
The reason for these problems
is that if
$f$ is not differentiable at $x$, then the gradient $\nabla f$ is discontinuous at $x$ and therefore $\nabla f(x)$ does not give any information about the behavior of $\nabla f$ in a neighborhood of $x$.

Not surprisingly, like in smooth optimization, the presence of constraints adds additional complexity, since constructing a descent sequence whose limit satisfies the constraints is (both theoretically and numerically) much more difficult than achieving this aim without the requirement of satisfying any restrictions.
~\\
\textbf{Linearly constrained nonsmooth optimization}. There exist various types of nonsmooth solvers like, e.g., the R-algorithm by \citet{Shor} or stochastic algorithms that try to approximate the subdifferential (e.g., by \citet{BLO}) 
or bundle algorithms which force a descent of the objective function by using local knowledge of the function. We will concentrate on the latter ones as they proved to be quite efficient.

One of the few publicly available bundle methods is the bundle-Newton method for nonsmooth, nonconvex unconstrained minimization by \citet{Luksan}. We sum up its key features:
It is the only method which we know of that uses second order information of the objective function, which results in faster convergence (in particular it was shown in \citet[p.~385,~Section~4]{Luksan} that the bundle-Newton method converges superlinearly for strongly convex functions).
Furthermore,
the
search direction is computed by solving a convex quadratic program (QP) (based on an SQP-approach in some sense)
and
it
uses a line search concept for deciding whether a serious step or a null step is performed.
Moreover,
its
implementation PNEW, which is described in \citet{LuksanPBUNandBNEW}, is written in FORTRAN.
Therefore, we can use the bundle-Newton method for solving linearly constrained nonsmooth optimization problems (as the linear constraints can just be inserted into the QP without any additional difficulties).

In general, every nonsmooth solver for unconstrained optimization can treat constrained problems via penalty functions. Nevertheless, choosing the penalty parameter well is a highly nontrivial task. Furthermore, if an application only allows the nonsmooth solver to perform a few steps (as, e.g., in \citet[\GenaueAngabeOne]{HannesPaperC}), we need to achieve a feasible descent within these steps.

\noindent\textbf{Nonlinearly constrained nonsmooth optimization}. Therefore, \citet[\GenaueAngabeANone]{HannesPaperB} give an extension of the bundle-Newton method to the constrained case in a very special way:
We use second order information of the constraint (cf.~(\ref{BundleSQP:OptProblem})).
Furthermore,
we use the SQP-approach of the bundle-Newton method for computing the
search direction for the constrained case and combine it with the idea
of quadratic constraint approximation, as it is used, e.g., in the
sequential quadratically constrained quadratic programming method by
\citet{Solodov} (this method is not a bundle method), in the hope to
obtain good feasible iterates, where we only accept strictly feasible
points as serious steps. Therefore, we have to solve a strictly
feasible convex QCQP for computing the search direction.
Using such a QCQP for computing the search direction yields a line
search condition for accepting infeasible points as trial points
(which is different to that in, e.g., \citet{MifflinNB}).
One of the most important properties of the convex QP (that is used to
determine the search direction) with respect to a bundle method is its
strong duality (e.g., for a meaningful termination criterion, for
global convergence,\dots) which is also true in the case of strictly
feasible convex QCQPs (cf.~\citet[\GenaueAngabeTwo]{HannesPaperB}).
Since there exist only a few solvers specialized in solving QCQPs (all
written in MATLAB or C, none in FORTRAN), the method is implemented in
MATLAB as well as in C.

For a detailed description of the presented issues we refer the reader to \citet[\GenaueAngabeThree]{HannesDissertation}.

The paper is organized as follows: In
Section~\ref{Paper:PresentationOfTheAlgorithm} we give a brief
description of the implemented variant of the second order bundle
algorithm. In Section~\ref{Paper:ReducedProblem} we discuss some
aspects that arise when using a convex QCQP for the computation of the
search direction problem like the reduction of its dimension and the
existence of a strictly feasible starting point for its
SOCP-reformulation.  Furthermore, we justify the approach for
determining the search direction by solving a QCQP numerically by
comparing the results of some well-known solvers for our search
direction problem. In Section~\ref{Paper:NumericalResults} we provide
numerical results for our second order bundle algorithm for some
examples of the Hock-Schittkowski collection by
\citet{Schittkowski,Schittkowski2}, for custom examples that arise in
the context of finding exclusion boxes for a quadratic CSP (constraint
satisfaction problem) in GloptLab by \citet{Domes} as well as for
higher dimensional piecewise quadratic examples, and finally we
compare these results to those of MPBNGC by \citet{MPBNGC} and SolvOpt
by \citet{SolvOptPublication} to emphasize the good performance of the
algorithm on constrained problems.

Throughout the paper we use the following notation:
We denote the non-negative real numbers by $\Rpos:=\lbrace
x\in\mathbb{R}:~x\geq0\rbrace$, and the space of all symmetric
$n\times n$-matrices by $\Sym{n}$. For $x\in\mathbb{R}^n$ we denote
the Euclidean norm of $x$ by $\lvert x\rvert$, for $1\leq i\leq j\leq
n$ we define the (MATLAB-like) colon operator $x_{i:j}:=(x_i,\dots
x_j)$, and for $A\in\mathrm{Sym(n)}$ we denote the spectral norm of
$A$ by $\lvert A\rvert$.

\section{Presentation of the algorithm}
\label{Paper:PresentationOfTheAlgorithm}
In the following section we give a brief exposition of our implemented
variant of the second order bundle algorithm whose theoretical
convergence properties are proved in \citet{HannesPaperB}.
For this purpose we assume that the functions
$f,F:\mathbb{R}^n\longrightarrow\mathbb{R}$ are locally Lipschitz
continuous, that $g_j\in\partial f(y_j)$ and $\hat{g}_j\in\partial
F(y_j)$, and
$G_j\in\subhessian{f}{y_j}$, $\hat{G}_j\in\subhessian{F}{y_j}$,
where the set $\subhessian{f}{x}\subseteq\Sym{n}$ of the substitutes
for the Hessian of $f$ at $x$ is defined by
\begin{equation*}
\subhessian{f}{x}
:=
\left\lbrace
\begin{array}{ll}
\lbrace G\rbrace & \textnormal{if the Hessian }G\textnormal{ of }f\textnormal{ at }x\textnormal{ exists}\\
\Sym{n}          & \textnormal{otherwise}\komma
\end{array}
\right.
\end{equation*}
i.e., we calculate elements of the sets $\subhessian{f}{y}$ and
$\subhessian{F}{y}$ (in the proof of convergence in \cite{HannesPaperB}
only approximations were required). We consider the nonsmooth
optimization problem (\ref{BundleSQP:OptProblem}) which has a single
nonsmooth constraint.
Then the second order bundle algorithm (described in Algorithm \ref{BundleSQPmitQCQP:Alg:GesamtAlgMitQCQP}) tries to solve optimization problem (\ref{BundleSQP:OptProblem}) according to the following scheme: After choosing a starting point $x_1\in\mathbb{R}^n$ and setting up a few positive definite matrices, we compute the
localized
approximation errors. Then we solve a convex QCQP to determine the search direction, where the intention of the usage of the quadratic constraints of the QCQP is to obtain preferably feasible points that yield a good descent.
\cbstartMIFFLIN
Therefore, we only use quadratic terms in the QCQP for the approximation of the constraint, but not for the approximation of the objective function (in contrast to \citet[\GenaueAngabeMIFFLIN]{HannesPaperB})
to balance the effort of solving the QCQP with the higher number of iterations caused by this simplification
(in Subsection \ref{PaperA:subsection:ReductionOfProblemSize} we will even discuss a further reduction of the size of the QCQP).
\cbendMIFFLIN
Now, after computing the aggregated data and the predicted descent as well as testing the termination criterion, we perform a line search (s.~Algorithm \ref{BundleSQP:AlgNB:LinesearchMitQCQP}) on the ray given by the search direction which yields a trial point $y_{k+1}$ that has the following property: Either $y_{k+1}$ is strictly feasible and the objective function achieves sufficient descent (serious step) or $y_{k+1}$ is strictly feasible and the model of the objective function changes sufficiently (null step with respect to the objective function) or $y_{k+1}$ is not strictly feasible and the model of the constraint changes sufficiently (null step with respect to the constraint). Afterwards we update the iteration point $x_{k+1}$ and the information which is stored in the bundle. Now, we repeat this procedure until the termination criterion is satisfied.
\begin{algorithm}
\label{BundleSQPmitQCQP:Alg:GesamtAlgMitQCQP}
\begin{enumerate}
\addtocounter{enumi}{-1}
\item\textit{Initialization:}

Choose the following parameters, which will not be changed during the algorithm:
\begin{longtable}{l|l|l}
\caption{Initial parameters}\\
\multicolumn{1}{c|}{\textnormal{General}}     & 
\multicolumn{1}{c|}{\textnormal{Default}}     & 
\multicolumn{1}{c }{\textnormal{Description}}  \\
\hline
\endfirsthead
\caption[]{Initial parameters (continued)}\\
\multicolumn{1}{c|}{\textnormal{General}}     & 
\multicolumn{1}{c|}{\textnormal{Default}}     & 
\multicolumn{1}{c }{\textnormal{Description}}  \\
\hline
\endhead
$x_1\in\mathbb{R}^n$       &                   & Strictly feasible initial point                    \\
$y_1=x_1$                  &                   & Initial trial point                                \\
$\varepsilon\geq0$         &                   & Final optimality tolerance                         \\
$M\geq2$                   & $M=n+3$           & Maximal bundle dimension                           \\
$t_0\in(0,1)$              & $t_0=0.001$       & Initial lower bound for step size                  \\
                           &                   & of serious step in line search                     \\
$\hat{t}_0\in(0,1)$        & $\hat{t}_0=0.001$ & Scaling parameter for $t_0$                        \\
$m_L\in(0,\tfrac{1}{2})$   & $m_L=0.01$        & Descent parameter for serious step in line search  \\
$m_R\in(m_L,1)$            & $m_R=0.5$         & Parameter for change of model of objective function\\
                           &                   & for short serious and null steps in line search    \\
$m_F\in(0,1)$              & $m_F=0.01$        & Parameter for change of model of constraint        \\
                           &                   & for short serious and null steps in line search    \\
$\zeta\in(0,\tfrac{1}{2})$ & $\zeta=0.01$      & Coefficient for interpolation in line search       \\
$\vartheta\geq1$           & $\vartheta=1$     & Exponent for interpolation in line search          \\
$C_S>0$                    & $C_S=10^{50}$     & Upper bound of the distance between $x_k$ and $y_k$\\
$C_G>0$                    & $C_G=10^{50}$     & Upper bound of the norm of the damped              \\
                           &                   & matrices $\lbrace\rho_jG_j\rbrace$ ($\lvert\rho_jG_j\rvert\leq C_G$)\\
$\hat{C}_G>0$              & $\hat{C}_G=C_G$   & Upper bound of the norm of the damped              \\
                           &                   & matrices $\lbrace\hat{\rho}_j\hat{G}_j\rbrace$ ($\lvert\hat{\rho}_j\hat{G}_j\rvert\leq \hat{C}_G$)\\
$\bar{\hat{C}}_G>0$        & $\bar{\hat{C}}_G=C_G$ & Upper bound of the norm of the matrices\\
                           &                       &
$\lbrace\bar{\hat{G}}_j^k\rbrace$ and $\lbrace\bar{\hat{G}}^k\rbrace$ ($\max{(\lvert\bar{\hat{G}}_j^k\rvert,\lvert\bar{\hat{G}}^k\rvert)}\leq\bar{\hat{C}}_G$)\\
$i_{\rho}\geq0$            & $i_{\rho}=3$          & Selection parameter for $\rho_{k+1}$
\\
$i_m\geq0$                 &                       & Matrix selection parameter
\\
$i_r\geq0$                 &                       & Bundle reset parameter
\\
$\gamma_1>0$               & $\gamma_1=1$      & Coefficient for locality measure for objective function\\
$\gamma_2>0$               & $\gamma_2=1$      & Coefficient for locality measure for constraint        \\
$\omega_1\geq1$            & $\omega_1=2$      & Exponent for locality measure for objective function   \\
$\omega_2\geq1$            & $\omega_2=2$      & Exponent for locality measure for constraint
\end{longtable}
Set the initial values of the data which gets changed during the algorithm:
\begin{align*}
i_n&=\hphantom{\lbrace}0\hphantom{\rbrace}\textnormal{ (\# subsequent null and short steps)}\\
i_s&=\hphantom{\lbrace}0\hphantom{\rbrace}\textnormal{ (\# subsequent serious steps)}       \\
J_1&=          \lbrace 1          \rbrace \textnormal{ (set of bundle indices)}             \punkt
\end{align*}
Compute the following information at the initial trial point
\begin{align*}
            f_p^1=f_1^1&=f(y_1)                                     
\\
            g_p^1=g_1^1&=g(y_1)\in\partial f(y_1)                   
\\
              G_p^1=G_1&=G(y_1)\in\subhessian{f}{y_1}               
\\
            F_p^1=F_1^1&=F(y_1)<0~~~\textnormal{(}y_1\textnormal{ is strictly feasible according to assumption)}
\\
\hat{g}_p^1=\hat{g}_1^1&=\hat{g}(y_1)\in\partial F(y_1)             
\\
  \hat{G}_p^1=\hat{G}_1&=\hat{G}(y_1)\in\subhessian{F}{y_1}    
\end{align*}
and set
\begin{align*}
\hat{s}_p^1=s_p^1=s_1^1&=0\textnormal{ (locality measure)}
\\
                \hat{\rho}_1=\rho_1&=1\textnormal{ (damping parameter)}                
\\
                     \bar{\kappa}^1&=1\textnormal{ (Lagrange multiplier for optimality condition)} 
\\
                                  k&=1\textnormal{ (iterator)}                          \punkt
\end{align*}
\item\textit{Determination of the matrices for the QCQP:}

\texttt{if} (step $k-1$ and $k-2$ were serious steps) $\wedge$ ($\lambda_{k-1}^{k-1}=1$ $\vee$ $\underset{\textnormal{bundle reset}}{\underbrace{i_s>i_r}}$)\\
\hphantom{aaa}$W=G_k+\bar{\kappa}^k\hat{G}_k$\\
\texttt{else}\\
\hphantom{aaa}$W=G_p^k+\bar{\kappa}^k\hat{G}_p^k$\\
\texttt{end}\\
\\
\texttt{if} $i_n\leq i_m$\\
\hphantom{aaa}$\widebar{W}_p^k=\textnormal{``positive definite modification of }W\textnormal{''}$\\
\texttt{else}\\
\hphantom{aaa}$\widebar{W}_p^k=\widebar{W}_p^{k-1}$\\
\texttt{end}\\
\cbstartDVI
\ifthenelse{\boolean{NurQCQPversion}}
{
Compute
\begin{equation}
\begin{split}
(\widebar{\hat{G}}^k,\widebar{\hat{G}}_j^k)=\textnormal{``positive definite modification of }(\hat{G}_p^k,\hat{G}_j)\textnormal{''}
~\forall j\in J_k
\punkt
\label{BundleSQPmitQCQP:Alg:ModifikationVonGhatbarjk}
\end{split}
\end{equation}
}
{
\texttt{if} $i_n<i_m+i_l$\\
\begin{equation}
\begin{split}
\widebar{\hat{G}}_j^k&=\textnormal{``positive definite modification of }
\hat{G}_j
\textnormal{''}\\
\widebar{\hat{G}}^k  &=\textnormal{``positive definite modification of }\hat{G}_p^k\textnormal{''}
\label{BundleSQPmitQCQP:Alg:ModifikationVonGhatbarjk}
\end{split}
\end{equation}
\texttt{else}\\
\hphantom{aaa}$\widebar{\hat{G}}_j^k=\zeroMatrix{n}$\\
\hphantom{aaa}$\widebar{\hat{G}}^k  =\zeroMatrix{n}$\\
\texttt{end}
}
\cbendDVI
\item\textit{Computation of the
localized
approximation errors:}
\begin{align*}
\alpha_j^k&:=\max{\big(\lvert f(x_k)-f_j^k\rvert,\gamma_1(s_j^k)^{\omega_1}\big)}
\komma\quad
\alpha_p^k:=\max\big(\lvert f(x_k) -f_p^k\rvert,\gamma_1(s_p^k)^{\omega_1}\big)
\\
A_j^k&:=\max{\big(\lvert F(x_k)-F_j^k\rvert,\gamma_2(s_j^k)^{\omega_2}\big)}
\komma\quad
A_p^k:=\max\big(\lvert F(x_k) -F_p^k\rvert,\gamma_2(\hat{s}_p^k)^{\omega_2}\big)\punkt
\end{align*}
\item\textit{Determination of the search direction:} Compute the solution $(d_k,\hat{v}_k)\in\mathbb{R}^{n+1}$ of the (convex) QCQP
\begin{equation}
\begin{split}
&\min_{d,\hat{v}}\hat{v}+\tfrac{1}{2}d^T\widebar{W}_p^kd\komma\\
&\textnormal{ s.t. }           -\alpha_j^k+d^Tg_j^k       \leq\hat{v}~~\hspace{82pt}\textnormal{for }j\in J_k\\
&\hphantom{\textnormal{ s.t. }}-\alpha_p^k+d^Tg_p^k       \leq\hat{v}~~\hspace{82pt}\textnormal{if }i_s\leq i_r\\
&\hphantom{\textnormal{ s.t. }}F(x_k)-A_j^k+d^T\hat{g}_j^k+\tfrac{1}{2}d^T\widebar{\hat{G}}_j^kd\leq0
~~\textnormal{for }j\in J_k\\
&\hphantom{\textnormal{ s.t. }}F(x_k)-A_p^k+d^T\hat{g}_p^k+\tfrac{1}{2}d^T\widebar{\hat{G}}^kd\leq0
~~\textnormal{if }i_s\leq i_r
\label{BundleSQPmitQCQP:Alg:QPTeilproblem}
\end{split}
\end{equation}
and its corresponding Lagrange multiplier $(\lambda^k,\lambda_p^k,\mu^k,\mu_p^k)\in\Rpos^{2(\lvert J_k\rvert+1)}$
and set
$
H_k:=
\big(\widebar{W}_p^k+\sum_{j\in J_k}\mu_j^k\widebar{\hat{G}}_j^k+\mu_p^k\widebar{\hat{G}}^k\big)
^{-\frac{1}{2}}
$ and
$
\bar{\kappa}^{k+1}
:=
\sum_{j\in J_k}\mu_j^k+\mu_p^k
$.\\
\texttt{if} $\bar{\kappa}^{k+1}>0$\\
\hphantom{aaa}$(\kappa_j^k,\kappa_p^k)=\tfrac{1}{\bar{\kappa}^{k+1}}(\mu_j^k,\mu_p^k)$\\
\texttt{else}\\
\hphantom{aaa}$(\kappa_j^k,\kappa_p^k)=0$\\
\texttt{end}\\
\texttt{if} $i_s>i_r$\\
\hphantom{aaa}$i_s=0$ (bundle reset)\\
\texttt{end}
\item\textit{Aggregation:} We set for the aggregation of information of the objective function
\begin{align*}
(\tilde{f}_p^k,\tilde{g}_p^k,G_p^{k+1},\tilde{s}_p^k)
&=
\sum\limits_{j\in J_k}\lambda_j^k(f_j^k,g_j^k,\rho_jG_j,s_j^k)+\lambda_p^k(f_p^k,g_p^k,G_p^k,s_p^k)
\\
\tilde{\alpha}_p^k&=\max\big(\vert f(x_k)-\tilde{f}_p^k\vert,\gamma_1(\tilde{s}_p^k)^{\omega_1}\big)
\end{align*}
and for the aggregation of information of the constraint
\begin{align*}
(\tilde{F}_p^k,\tilde{\hat{g}}_p^k,\hat{G}_p^{k+1},\tilde{\hat{s}}_p^k)
&=
\sum\limits_{j\in J_k}\kappa_j^k(F_j^k,\hat{g}_j^k,\hat{\rho}_j\hat{G}_j,s_j^k)
+\kappa_p^k(F_p^k,\hat{g}_p^k,\hat{G}_p^k,\hat{s}_p^k)
\\
      \tilde{A}_p^k&=\max\big(\vert F(x_k)-\tilde{F}_p^k\vert,\gamma_2(\tilde{\hat{s}}_p^k)^{\omega_2}\big)
\end{align*}
and we set
\begin{align*}
v_k
&=
-d_k^T\widebar{W}_p^kd_k
-\tfrac{1}{2}d_k^T\big(\sum_{j\in J_k}\mu_j^k\widebar{\hat{G}}_j^k+\mu_p^k\widebar{\hat{G}}^k\big)d_k
-\tilde{\alpha}_p^k-\bar{\kappa}^{k+1}\tilde{A}_p^k
-\bar{\kappa}^{k+1}\big(-F(x_k)\big)
\\
w_k
&=
\tfrac{1}{2}\vert H_k
(
\tilde{g}_p^k+\bar{\kappa}^{k+1}\tilde{\hat{g}}_p^k
)
\vert^2
+\tilde{\alpha}_p^k+\bar{\kappa}^{k+1}\tilde{A}_p^k
+\bar{\kappa}^{k+1}\big(-F(x_k)\big)
\punkt
\end{align*}
\item\textit{Termination criterion:}

\texttt{if} $w_k\leq\varepsilon$\\
\hphantom{aaa}\texttt{stop}\\
\texttt{end}

\item\textit{Line search:} We compute step sizes $0\leq t_L^k\leq t_R^k\leq1$ and $t_0^k\in(0,t_0]$ by using the line search described in Algorithm \ref{BundleSQP:AlgNB:LinesearchMitQCQP}
and we set
\begin{align*}
      x_{k+1}&=x_k+t_L^kd_k~~~\textnormal{(is created strictly feasible by the line search)}
\\
      y_{k+1}&=x_k+t_R^kd_k                              
\\
      f_{k+1}&=f(y_{k+1})                                
\komma\quad
      g_{k+1}=g(y_{k+1})\in\partial f(y_{k+1})          
\komma\quad
      G_{k+1}=G(y_{k+1})\in\subhessian{f}{y_{k+1}}      
\\
      F_{k+1}&=F(y_{k+1})                                
\komma\quad
\hat{g}_{k+1}=\hat{g}(y_{k+1})\in\partial F(y_{k+1})    
\komma\quad
\hat{G}_{k+1}=\hat{G}(y_{k+1})\in\subhessian{F}{y_{k+1}}\punkt
\end{align*}
\item\textit{Update:}

\texttt{if} $i_n\leq i_{\rho}$\\
\hphantom{aaa}$\rho_{k+1}=\min(1,\tfrac{C_G}{\vert G_{k+1}\vert})$\\
\texttt{else}\\
\hphantom{aaa}$\rho_{k+1}=0$\\
\texttt{end}\\
$\hat{\rho}_{k+1}=\min(1,\tfrac{\hat{C}_G}{\vert\hat{G}_{k+1}\vert})$\\

\texttt{if} $t_L^k\geq t_0^k$ (serious step)\\
\hphantom{aaa}$i_n=0$            \\
\hphantom{aaa}$i_s=i_s+1$        \\
\texttt{else} (no serious step, i.e.~null or short step)
\\
\hphantom{aaa}$i_n=i_n+1$\\
\texttt{end}\\
\\
Compute the updates of the locality measure
\begin{align*}
    s_j^{k+1}&=s_j^k+\vert x_{k+1}-x_k\vert~~~\textnormal{for }j\in J_k
\\
s_{k+1}^{k+1}&=\vert x_{k+1}-y_{k+1}\vert
\\
    s_p^{k+1}&=\tilde{s}_p^k+\vert x_{k+1}-x_k\vert  
\\
    \hat{s}_p^{k+1}&=\tilde{\hat{s}}_p^k+\vert x_{k+1}-x_k\vert\punkt
\end{align*}
Compute the updates for the objective function approximation
\begin{align*}
f_j^{k+1}
&=
f_j^k+g_j^{k\,T}(x_{k+1}-x_k)+\tfrac{1}{2}\rho_j(x_{k+1}-x_k)^TG_j(x_{k+1}-x_k)~~~\textnormal{for }j\in J_k
\\
f_{k+1}^{k+1}
&=
f_{k+1}+g_{k+1}^T(x_{k+1}-y_{k+1})+\tfrac{1}{2}\rho_{k+1}(x_{k+1}-y_{k+1})^TG_{k+1}(x_{k+1}-y_{k+1})
\\
f_p^{k+1}
&=
\tilde{f}_p^k+\tilde{g}_p^{k\,T}(x_{k+1}-x_k)+\tfrac{1}{2}(x_{k+1}-x_k)^TG_p^{k+1}(x_{k+1}-x_k)
\end{align*}
and for the constraint
\begin{align*}
F_j^{k+1}
&=
F_j^k+\hat{g}_j^{k\,T}(x_{k+1}-x_k)+\tfrac{1}{2}\hat{\rho}_j(x_{k+1}-x_k)^T\hat{G}_j(x_{k+1}-x_k)
~~~\textnormal{for }j\in J_k
\\
F_{k+1}^{k+1}
&=
F_{k+1}+\hat{g}_{k+1}^T(x_{k+1}-y_{k+1})+\tfrac{1}{2}\hat{\rho}_{k+1}(x_{k+1}-y_{k+1})^T\hat{G}_{k+1}(x_{k+1}-y_{k+1})
\\
F_p^{k+1}
&=
\tilde{F}_p^k+\tilde{\hat{g}}_p^{k\,T}(x_{k+1}-x_k)+\tfrac{1}{2}(x_{k+1}-x_k)^T\hat{G}_p^{k+1}(x_{k+1}-x_k)
\punkt
\end{align*}
Compute the updates for the subgradient of the objective function approximation
\begin{align*}
g_j^{k+1}
&=
g_j^k+\rho_jG_j(x_{k+1}-x_k)~~~\textnormal{for }j\in J_k
\\
g_{k+1}^{k+1}
&=
g_{k+1}+\rho_{k+1}G_{k+1}(x_{k+1}-y_{k+1})
\\
g_p^{k+1}
&=
\tilde{g}_p^k+G_p^{k+1}(x_{k+1}-x_k)
\end{align*}
and for the constraint
\begin{align*}
\hat{g}_j^{k+1}
&=
\hat{g}_j^k+\hat{\rho}_j\hat{G}_j(x_{k+1}-x_k)
~~~\textnormal{for }j\in J_k
\\
\hat{g}_{k+1}^{k+1}
&=
\hat{g}_{k+1}+\hat{\rho}_{k+1}\hat{G}_{k+1}(x_{k+1}-y_{k+1})
\\
\hat{g}_p^{k+1}
&=
\tilde{\hat{g}}_p^k+\hat{G}_p^{k+1}(x_{k+1}-x_k)
\punkt
\end{align*}
Choose $J_{k+1}\subseteq\lbrace k-M+2,\dots,k+1\rbrace\cap\lbrace1,2,\dots\rbrace$ with $k+1\in J_{k+1}$.\\
$k=k+1$\\
Go to 1
\end{enumerate}
\end{algorithm}
We extend the line search of the bundle-Newton method for nonsmooth unconstrained minimization to the constrained case in the line search described in Algorithm \ref{BundleSQP:AlgNB:LinesearchMitQCQP}.
Before formulating the line search in detail, we give a brief overview of its functionality:

Starting with the step size $t=1$, we check if the point $x_k+td_k$ is strictly feasible.
If so and if additionally the objective function decreases sufficiently in this point and $t$ is not too small, then we take $x_k+td_k$ as new iteration point in Algorithm \ref{BundleSQPmitQCQP:Alg:GesamtAlgMitQCQP} (serious step).
Otherwise, if the point $x_k+td_k$ is strictly feasible and the model of the objective function changes sufficiently, we take $x_k+td_k$ as new trial point (short/null step with respect to the objective function).
If $x_k+td_k$ is not strictly feasible, but the model of the constraint changes sufficiently (in particular here the quadratic approximation of the constraint comes into play), we take $x_k+td_k$ as new trial point (short/null step with respect to the constraint).
After choosing a new step size $t\in[0,1]$ by interpolation, we iterate this procedure.
\begin{algorithm}
\label{BundleSQP:AlgNB:LinesearchMitQCQP}
\begin{enumerate}
\addtocounter{enumi}{-1}
\item\textit{Initialization:} Choose $\zeta\in(0,\tfrac{1}{2})$ as well as $\vartheta\geq1$ and set $t_L=0$ as well as $t=t_U=1$.
\item\textit{Modification of either $t_L$ or $t_U$:}
\begin{align*}
&\texttt{if }F(x_k+td_k)<0                                            \nonumber\\
&\hphantom{iii}\texttt{if }f(x_k+td_k)\leq f(x_k)+m_Lv_k\cdot t       \nonumber\\
&\hphantom{iiiiii}t_L=t                                                     \nonumber\\
&\hphantom{iii}\texttt{else if }f(x_k+td_k)>f(x_k)+m_Lv_k\cdot t       \nonumber\\
&\hphantom{iiiiii}t_U=t                                                     \nonumber\\
&\hphantom{iii}\texttt{end}                                                 \nonumber\\
\end{align*}
\begin{align*}
&\texttt{else if }F(x_k+td_k)\geq0                                     \nonumber\\
&\hphantom{iii}t_U=t                                                        \nonumber\\
&\hphantom{iii}t_0=\hat{t}_0t_U
\\
&\texttt{end}                                                               \nonumber\\
&\texttt{if }t_L\geq t_0                                               \nonumber\\
&\hphantom{iii}t_R=t_L                                                      \nonumber\\
&\hphantom{iii}\texttt{return}\textnormal{ (serious step)}                  \nonumber\\
&\texttt{end}                                                               \nonumber
\end{align*}
\item\textit{Decision of return:}
\begin{align}
&\hphantom{\texttt{10 }}\texttt{if }F(x_k+td_k)<0
\nonumber\\
&\hphantom{\texttt{10 }}\hphantom{iii}
g=g(x_k+td_k)\in\partial f(x_k+td_k)
\komma\quad
G=G(x_k+td_k)\in\subhessian{f}{x_k+td_k}\nonumber\\
&\hphantom{\texttt{10 }}\hphantom{iii}
\rho=
\left\lbrace
\begin{array}{ll}
\min(1,\tfrac{C_G}{\lvert G\rvert}) & \textnormal{for }i_n\leq 3\\
0                                   & \textnormal{else}
\end{array}
\right.                                                                                                              
\nonumber\\
&\hphantom{\texttt{10 }}\hphantom{iii}
f=f(x_k+td_k)+(t_L-t)g^Td_k+\tfrac{1}{2}\rho(t_L-t)^2d_k^TGd_k
\nonumber
\\
&\hphantom{\texttt{10 }}\hphantom{iii}
\beta=
\max(\lvert f(x_k+t_Ld_k)-f\rvert,\gamma_1\lvert t_L-t\rvert^{\omega_1}\lvert d_k\rvert^{\omega_1})
\nonumber
\\
&\hphantom{\texttt{10 }}\hphantom{iii}\texttt{if }
-\beta+d_k^T\big(g+\rho(t_L-t)Gd_k\big)\geq m_Rv_k\texttt{ and }(t-t_L)\lvert d_k\rvert\leq C_S
\nonumber\\
&\hphantom{\texttt{10 }}\hphantom{iiiiii}t_R=t
\nonumber\\
&\hphantom{\texttt{10 }}\hphantom{iiiiii}\texttt{return}\textnormal{ (short/null step: change of model of the objective function)}
\nonumber\\
&\hphantom{\texttt{10 }}\hphantom{iii}\texttt{end}
\nonumber\\
&\hphantom{\texttt{10 }}\texttt{else if }F(x_k+td_k)\geq0
\nonumber\\
&\hphantom{\texttt{10 }}\hphantom{iii}
\hat{g}=\hat{g}(x_k+td_k)\in\partial F(x_k+td_k)
\komma\quad
\hat{G}=\hat{G}(x_k+td_k)\in\subhessian{F}{x_k+td_k}\nonumber\\
&\hphantom{\texttt{10 }}\hphantom{iii}
\hat{\rho}=
\min(1,\tfrac{\hat{C}_G}{\lvert\hat{G}\rvert})
\nonumber\\
&\hphantom{\texttt{10 }}\hphantom{iii}
F=F(x_k+td_k)+(t_L-t)\hat{g}^Td_k+\tfrac{1}{2}\rho(t_L-t)^2d_k^T\hat{G}d_k
\nonumber
\\
&\hphantom{\texttt{10 }}\hphantom{iii}
\hat{\beta}=
\max(\lvert F(x_k+t_Ld_k)-F\rvert,\gamma_2\lvert t_L-t\rvert^{\omega_2}\lvert d_k\rvert^{\omega_2})
\nonumber
\\
&\hphantom{\texttt{10 }}\hphantom{iii}
\widebar{\hat{G}}=\textnormal{``positive definite modification of }\hat{G}\textnormal{''}
\label{Luk:LineareLS3:NEWinfeasibleTrialPointPositiveDefiniteModification}
\\
&\hphantom{\texttt{10 }}\hphantom{iii}\texttt{if }
F(x_k+t_Ld_k)-\hat{\beta}+d_k^T\big(\hat{g}+\hat{\rho}(t_L-t)\hat{G}d_k\big)\geq m_F\cdot(-\tfrac{1}{2}d_k^T\widebar{\hat{G}}d_k)\nonumber\\
&\hspace{5cm}\texttt{ and }(t-t_L)\lvert d_k\rvert\leq C_S
\label{Luk:LineareLS3:NEWinfeasibleTrialPoint}
\\
&\hphantom{\texttt{10 }}\hphantom{iiiiii}t_R=t
\nonumber\\
&\hphantom{\texttt{10 }}\hphantom{iiiiii}\texttt{return}\textnormal{ (short/null step: change of model of the constraint)}
\nonumber\\
&\hphantom{\texttt{10 }}\hphantom{iii}\texttt{end}
\nonumber\\
&\hphantom{\texttt{10 }}\texttt{end}
\nonumber
\end{align}
\item\textit{Interpolation:} Choose
$
t\in[t_L+\zeta(t_U-t_L)^{\vartheta},t_U-\zeta(t_U-t_L)^{\vartheta}]
$.
\item\textit{Loop:} \texttt{Go to 1}
\end{enumerate}
\end{algorithm}
\begin{remark} 
Similar to the line search in the bundle-Newton method for nonsmooth unconstrained minimization by \citet{Luksan}, we want to choose a new point in the interval $[t_L+\zeta(t_U-t_L)^{\vartheta},t_U-\zeta(t_U-t_L)^{\vartheta}]$ by interpolation.
For this purpose, we set up a polynomial $p$ passing through $\big(t_L,f(x_k+t_Ld_k)\big)$ and $\big(t_U,f(x_k+t_Ud_k)\big)$ as well as a polynomial $q$ passing through $\big(t_L,F(x_k+t_Ld_k)\big)$ and $\big(t_U,F(x_k+t_Ud_k)\big)$. Now we minimize $p$ subject to the constraint $q(t)\leq0$ on $[t_L+\zeta(t_U-t_L)^{\vartheta},t_U-\zeta(t_U-t_L)^{\vartheta}]$
and we use a solution $\hat{t}$
as the new point. The degree of the polynomial should be chosen in a way that
determining $\hat{t}$
is easy (e.g., if we choose $p$ and $q$ as quadratic polynomials, then
determining $\hat{t}$
consists of solving a one-dimensional linear equation, a one-dimensional quadratic equation and a few case distinctions).
\end{remark} 

\section{The reduced problem}
\label{Paper:ReducedProblem}
In this section we present some issues that arise when using a convex QCQP for the computation of the search direction problem like the reduction of its dimension.
Moreover, we give a numerical justification of the approach of determining the search direction by solving a QCQP by comparing the results of some well-known solvers for our search direction problem.
\subsection{Reduction of problem size}
\label{PaperA:subsection:ReductionOfProblemSize}
We want to reduce the problem size of the QCQP (\ref{BundleSQPmitQCQP:Alg:QPTeilproblem}). For this purpose we choose
$\widebar{\hat{G}}^k$ as a positive definite modification of $\hat{G}_p^k$ and
$\widebar{\hat{G}}_j^k:=\widebar{\hat{G}}^k$ for all $j\in J_k$,
i.e.~we choose all matrices for the constraint approximation equal to a positive definite modification of an aggregated Hessian of the constraint (i.e.~similar to the choice of $\widebar{W}_p^k$ in the bundle-Newton method for nonsmooth unconstrained minimization by \citet{Luksan}). For the implementation, we will extract linear constraints $Bx\leq b$ with $B\in\mathbb{R}^{\bar{m}\times n}$ and $b\in\mathbb{R}^{\bar{m}}$ that may occur in the single nonsmooth function $F:\mathbb{R}^n\longrightarrow\mathbb{R}$ (via a $\max$-function of the rows $B_{i:}x-b_i\leq0$ for all $i=1,\dots,\bar{m}$) in the nonsmooth constrained optimization problem (\ref{BundleSQP:OptProblem}) and put them directly into the search direction problem (this is the usual way of handling linear constraints in bundle methods). For easiness of exposition, we drop the p-constraints. These facts altogether yield the (convex) QCQP
\begin{equation}
\begin{split}
&\min_{d,\hat{v}}\hat{v}+\frac{1}{2}d^T\widebar{W}_p^kd\\
&\textnormal{ s.t. }           -\alpha_j^k+d^Tg_j^k       \leq\hat{v}~~\hspace{81.5pt}\textnormal{for }j\in J_k\\
&\hphantom{\textnormal{ s.t. }}F(x_k)-A_j^k+d^T\hat{g}_j^k+\tfrac{1}{2}d^T\widebar{\hat{G}}^kd\leq0
~~\textnormal{for }j\in J_k\\
&\hphantom{\textnormal{ s.t. }}B_{i:}(x_k+d)\leq b_i~~\hspace{88pt}\textnormal{for }i=1,\dots,\bar{m}
\punkt
\label{BundleSQPmitQCQP:Alg:QPTeilproblemReduced}
\end{split}
\end{equation}
Furthermore, we consider the following
modification
of the QCQP (\ref{BundleSQPmitQCQP:Alg:QPTeilproblemReduced})
\begin{equation}
\begin{split}
&\min_{d,\hat{v},\hat{u}}\hat{v}+\frac{1}{2}d^T\widebar{W}_p^kd\\
&\textnormal{ s.t. }           -\alpha_j^k+d^Tg_j^k               \leq\hat{v}~~\hspace{49.3pt}\textnormal{for }j\in J_k\\
&\hphantom{\textnormal{ s.t. }}F(x_k)-A_j^k+d^T\hat{g}_j^k+\hat{u}\leq0
~~\textnormal{for }j\in J_k\\
&\hphantom{\textnormal{ s.t. }}\tfrac{1}{2}d^T\widebar{\hat{G}}^kd    \leq\hat{u}\\
&\hphantom{\textnormal{ s.t. }}B_{i:}(x_k+d)                      \leq b_i
~~\hspace{56pt}\textnormal{for }i=1,\dots,\bar{m}\komma
\label{SOCP:Luksan:Alg:QCQPTeilproblem2QsocpKapitel}
\end{split}
\end{equation}
which is a (convex) QCQP with only one quadratic constraint.
\begin{remark}
We expect that the reduced QCQP (\ref{SOCP:Luksan:Alg:QCQPTeilproblem2QsocpKapitel}) should be solved much faster than the QCQP (\ref{BundleSQPmitQCQP:Alg:QPTeilproblemReduced}) because
of the following reasons:

An interior point method for solving QPs/QCQPs solves a linear system (called the KKT-system) at each iteration which is the most time consuming operation, i.e.~the bigger the KKT-system is, the longer the interior point method will need to solve the problem.

If we solved a QP to determine the search direction (we do not do this because of \citet[\GenaueAngabeFour]{HannesPaperB}), we would obtain $\lvert J_k\rvert+1$ linear constraints for approximating $F$ which increases the size of the KKT-system by $\lvert J_k\rvert+1$ rows compared to the unconstrained case (i.e.~without $F$).

If we solve the QCQP (\ref{BundleSQPmitQCQP:Alg:QPTeilproblemReduced}) to determine the search direction, we will obtain --- in addition to the $\lvert J_k\rvert+1$ rows which are due to the linear terms --- $\lvert J_k\rvert+1$ many $n\times n$-blocks (i.e.~$(\lvert J_k\rvert+1)n$ rows) which are due to the $\lvert J_k\rvert+1$ quadratic terms. Since $J_k$ is bounded by the maximal bundle dimension $M$ and if we choose, e.g., $M=n+3$ (this is the recommended default value for $M$ in the bundle-Newton method by \citet{Luksan} for nonsmooth unconstrained minimization), then the KKT-system can become
very
big
even for low dimensions.

If we solve the reduced QCQP (\ref{SOCP:Luksan:Alg:QCQPTeilproblem2QsocpKapitel}) to determine the search direction, we will obtain --- in addition to the $\lvert J_k\rvert+1$ rows which are due to the linear terms --- only one $n\times n$-block (i.e.~$n$ rows) since we only have one quadratic term. Therefore, if $n$ is not too big, we expect that solving the reduced QCQP should not take significantly more time than solving the corresponding QP at least for a good interior point method and this turns out to be true indeed
(cf.~the comparisons in Subsection \ref{subsection:Comparisons}).

So the big advantage of the reduced QCQP (\ref{SOCP:Luksan:Alg:QCQPTeilproblem2QsocpKapitel}) is that it has a size similar to that of the corresponding QP (i.e.~its size is much smaller than that of the QCQP (\ref{BundleSQPmitQCQP:Alg:QPTeilproblemReduced})), but it still uses quadratic information to deal with the nonlinearity of $F$.

Furthermore,
we
do not need to compute a positive definite modification $\widebar{\hat{G}}_j^k$ of
$\hat{G}_j$
in (\ref{BundleSQPmitQCQP:Alg:ModifikationVonGhatbarjk}), and
we
can replace the model change condition in (\ref{Luk:LineareLS3:NEWinfeasibleTrialPoint}) by
\begin{equation*}
F(x_k+t_Ld_k)-\hat{\beta}+d_k^T\big(\hat{g}+\hat{\rho}(t_L-t)\hat{G}d_k\big)\geq m_F\cdot(-\hat{u}_k)
\end{equation*}
and therefore we do not need to compute a positive definite modification $\widebar{\hat{G}}$ of $\hat{G}$ in (\ref{Luk:LineareLS3:NEWinfeasibleTrialPointPositiveDefiniteModification}).
\end{remark}
\subsection{Overview of the QCQP-solvers}
The most time-consuming part of the bundle-Newton method for nonsmooth unconstrained minimization by \citet{Luksan} is solving a (convex) QP. This QP is solved by the FORTRAN solver PLQDF1 described in \citet{LuksanQP} which exploits the special structure of the QP.
Analogously, the most time-consuming part of Algorithm \ref{BundleSQPmitQCQP:Alg:GesamtAlgMitQCQP} is solving the (convex) QCQP (\ref{BundleSQPmitQCQP:Alg:QPTeilproblem}).

For solving the QCQP (\ref{BundleSQPmitQCQP:Alg:QPTeilproblem}),
our implementation of
Algorithm \ref{BundleSQPmitQCQP:Alg:GesamtAlgMitQCQP} can use
MOSEK by \citet{AndersenPaper,AndersenMosek} (which is written in C and available as commercial software resp.~as a trial version without any limitations of the problem size that may be used by an academic institution for 90 days)
or
IPOPT by \citet{WaechterPaper,WaechterIpopt} (which is written in C++ and
freely available),
where the ordering represents the performance of the solvers according to the tests in \citet{Mittelmann1}.

For solving the SOCP-reformulation of the QCQP (\ref{BundleSQPmitQCQP:Alg:QPTeilproblem})
(cf.~\citet[\GenaueAngabeNEWone]{HannesDissertation} for details),
our implementation of
Algorithm \ref{BundleSQPmitQCQP:Alg:GesamtAlgMitQCQP} can use
MOSEK,
SEDUMI by \citet{Sturm,Polik} (which is written in MATLAB and freely available)
SDPT3 by \citet{SDPT3} (which is written in MATLAB and freely available), or
\texttt{socp} by \citet{BoydSOCP1995} (which is written in C and freely available).
Again, the ordering represents the performance of the solvers according to the tests in \citet{Mittelmann2}, except for \texttt{socp} which was not tested there.

The comparisons in \citet{Mittelmann1,Mittelmann2} coincide with our own observations
(cf.~Subsection \ref{subsection:Comparisons}).

\subsection{Comparison of the QCQP-solvers}
\label{subsection:Comparisons}
All tests were performed on an Intel Pentium IV with 3 GHz and 1 GB RAM running Microsoft Windows XP and MATLAB R2010a.

We are comparing the time for solving 50 randomly generated problems of the following types
\begin{align*}
\textnormal{L(inear)}&:=\textnormal{``QP obtained by setting }
\widebar{\hat{G}}_j^k=\widebar{\hat{G}}^k=\zeroMatrix{n}\textnormal{ in QCQP }(\ref{BundleSQPmitQCQP:Alg:QPTeilproblemReduced})\textnormal{''}
\\
\textnormal{D(ifferent)}&:=\textnormal{``QCQP }(\ref{BundleSQPmitQCQP:Alg:QPTeilproblemReduced})\textnormal{''}
\\
\textnormal{E(qual)}&:=\textnormal{``QCQP }(\ref{BundleSQPmitQCQP:Alg:QPTeilproblemReduced})\textnormal{ with }\widebar{\hat{G}}_j^k=\widebar{\hat{G}}^k\textnormal{''}\\
\textnormal{R(educed)}&:=\textnormal{``Reduced QCQP }(\ref{SOCP:Luksan:Alg:QCQPTeilproblem2QsocpKapitel})\textnormal{''}
\komma
\end{align*}
where we set $m:=\lvert J_k\rvert$ and we choose $\bar{m}=0$.

For obtaining a first insight, how long the computation of the search
direction will take, we compare the plots (based on the data from
Table~2 in appendix B of
\ifthenelse{\boolean{AppendixOutsourcing}}
{\citet{HannesPaperAOutsourcing})}
{Appendix \ref{section:ResultTables})}
of the median solving times (in milliseconds) for
the
MOSEK QCQP-solver ($\square$),
the
MOSEK SOCP-solver ($\diamondsuit$),
SEDUMI ($\nabla$), and
SDPT3 ($\triangle$),
where we use the
symbols to distinguish the results of the different solvers (since the
only purpose of this subsection is to obtain a rough estimation of the
solving times of the different types of search direction problems, we
only tested these solvers here because the MATLAB tools CVX by
\citet{BoydCVX} resp.~YALMIP by \citet{LoefbergYALMIP} offer an
excellent interface for easily generation of the input data of the
different search direction problems for these different solvers; the
performance of \texttt{socp} resp.~IPOPT is discussed in Remark
\ref{remark:socpBadMOSEKrobust} within the framework of using one of
these two algorithms as the (QC)QP-solver in Algorithm
\ref{BundleSQPmitQCQP:Alg:GesamtAlgMitQCQP}).
~\\
~\\
\begin{figure}[htb]
  \begin{minipage}{.5\linewidth}
    \begin{center}
      \captionsetup{type=figure}
      \includegraphics[width=6cm]{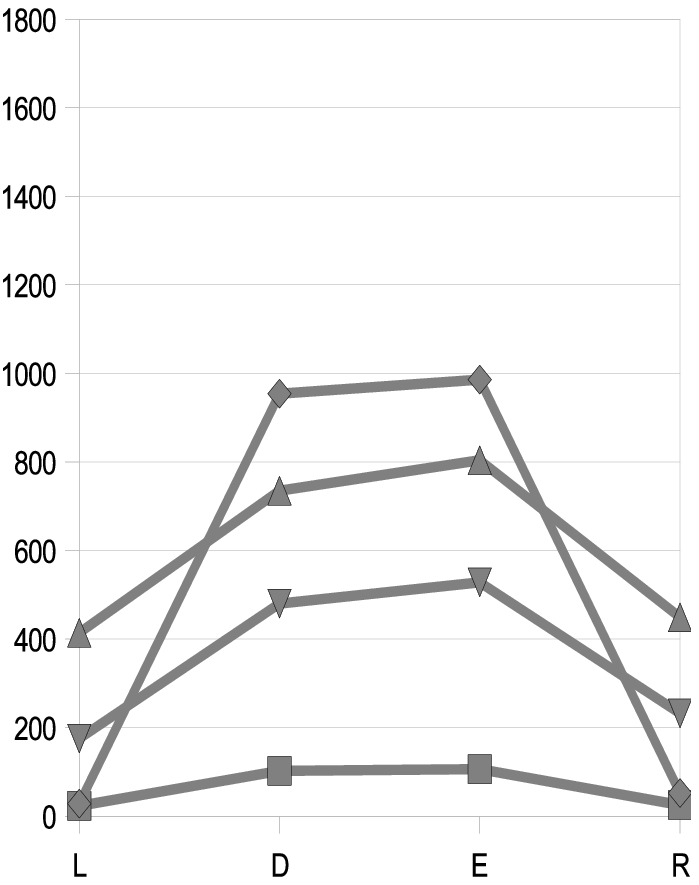}
      \captionof{figure}{Median solving time for $n=50$ and $m=25$}
      \label{Figure:ComparisonSDproblemAllSolvers[n=50][m=25]}
    \end{center}
  \end{minipage}
  \begin{minipage}{.5\linewidth}
    \begin{center}
      \captionsetup{type=figure}
      \includegraphics[width=6cm]{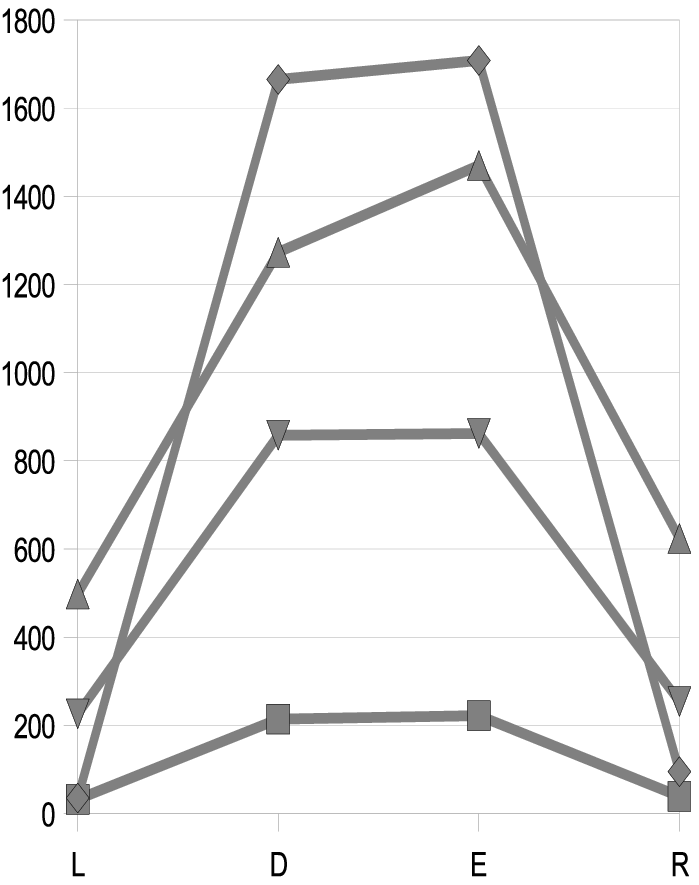}
      \captionof{figure}{Median solving time for $n=50$ and $m=50$}
      \label{Figure:ComparisonSDproblemAllSolvers[n=50][m=50]}
    \end{center}
  \end{minipage}
\end{figure}
\begin{figure}[htb]
  \begin{minipage}{.5\linewidth}
    \begin{center}
      \captionsetup{type=figure}
      \includegraphics[width=6cm]{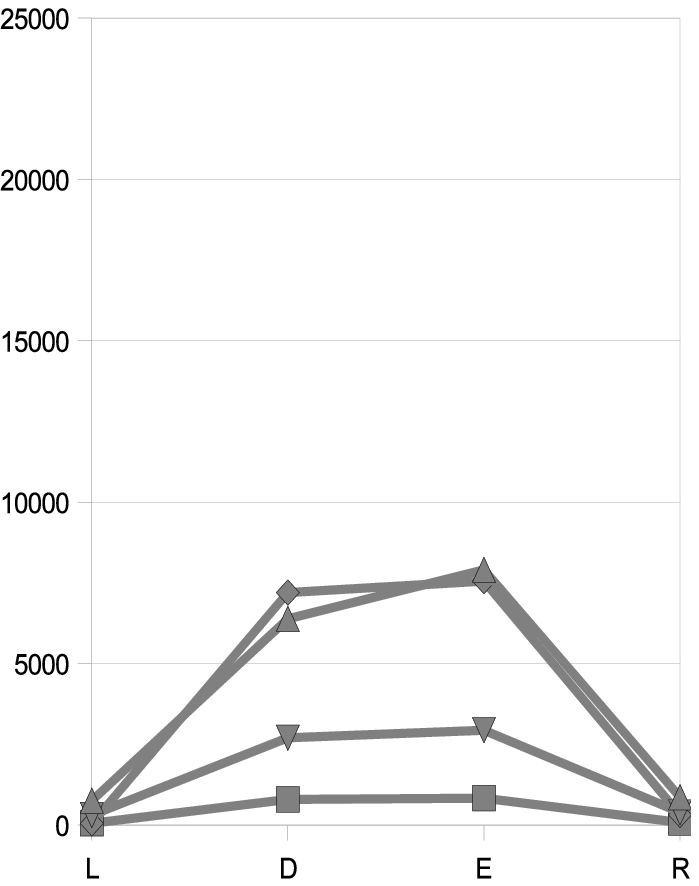}
      \captionof{figure}{Median solving time for $n=100$ and $m=50$}
      \label{Figure:ComparisonSDproblemAllSolvers[n=100][m=50]}
    \end{center}
  \end{minipage}
  \begin{minipage}{.5\linewidth}
    \begin{center}
      \captionsetup{type=figure}
      \includegraphics[width=6cm]{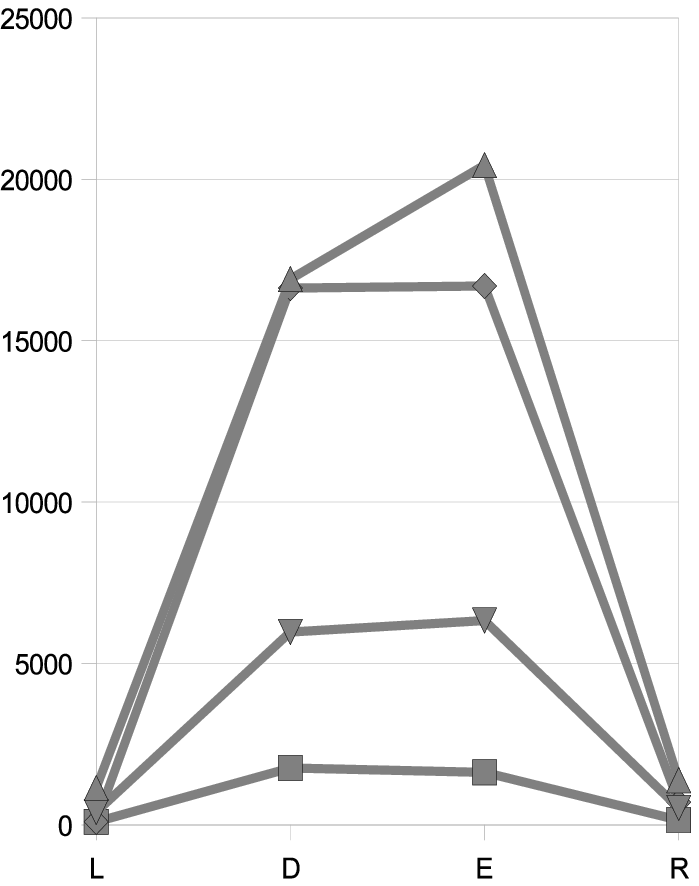}
      \captionof{figure}{Median solving time for $n=100$ and $m=100$}
      \label{Figure:ComparisonSDproblemAllSolvers[n=100][m=100]}
    \end{center}
  \end{minipage}
\end{figure}
In Figures~\ref{Figure:ComparisonSDproblemAllSolvers[n=50][m=25]}
to~\ref{Figure:ComparisonSDproblemAllSolvers[n=100][m=100]} we plot
the median of the solving times for various problem sizes. Here we see
that L and R are significantly faster than D and E. To analyze the
difference between the first two algorithms we magnify the results of
L (dashed line)
and
R (solid line)
and plot the result in Figure \ref{Figure:ComparisonSDproblemAllSolversMagnification}.
\begin{figure}[htb]
  \begin{center}
    \captionsetup{type=figure}
    \includegraphics[width=7cm]{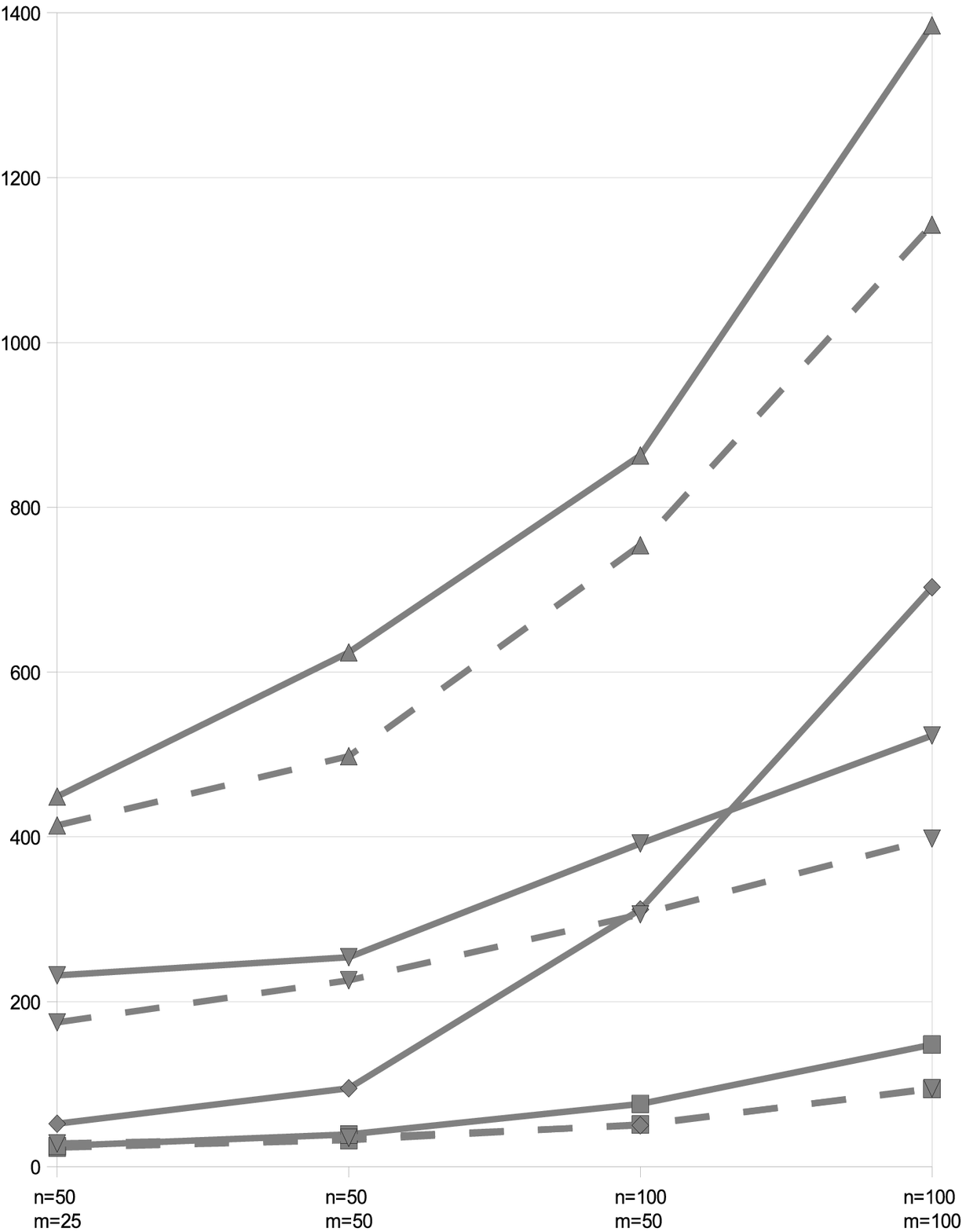}
    \captionof{figure}{Magnification of the median solving time for L
      and R}
    \label{Figure:ComparisonSDproblemAllSolversMagnification}
  \end{center}
\end{figure}
\begin{remark}
\label{Remark:ComparisonsOfMOSEK:QCQPandSOCP}
Although \citet[p.~131,~Section~7.2~and~7.2.1]{AndersenMosek} recommends to rather use the MOSEK SOCP-solver than the MOSEK QCQP-solver for solving convex QCQPs, this does not coincide with the above results in which the MOSEK QCQP-solver has a significantly better performance than the MOSEK SOCP-solver for solving a QCQP of our shape.
\end{remark}
The results from
Figures \ref{Figure:ComparisonSDproblemAllSolvers[n=50][m=25]}--\ref{Figure:ComparisonSDproblemAllSolversMagnification}
suggest that we will only test the MOSEK QCQP-solver on the reduced QCQP (\ref{SOCP:Luksan:Alg:QCQPTeilproblem2QsocpKapitel}) in higher dimensions as this is the only combination that does not significantly exceed the shortest duration for solving the corresponding QP (which is always achieved by the MOSEK QP-solver).
Therefore, we plot in Figure \ref{Figure:ComparisonSDproblemMedMinMax} (based on the data from Table 3 in appendix B of
\ifthenelse{\boolean{AppendixOutsourcing}}
{\citet{HannesPaperAOutsourcing})}
{Appendix \ref{section:ResultTables})}
the
minimal \& maximal (lower and upper end of the vertical line)
and the
median (horizontal line)
solving times (in milliseconds) obtained by MOSEK for
L
(black)
and
R
(grey)
(from $n=m=400$ on, our computer started to swap and, consequently, we did not test higher dimensional problems).
\begin{figure}[htb]
  \begin{center}
    \captionsetup{type=figure}
    \includegraphics[width=14.8cm]{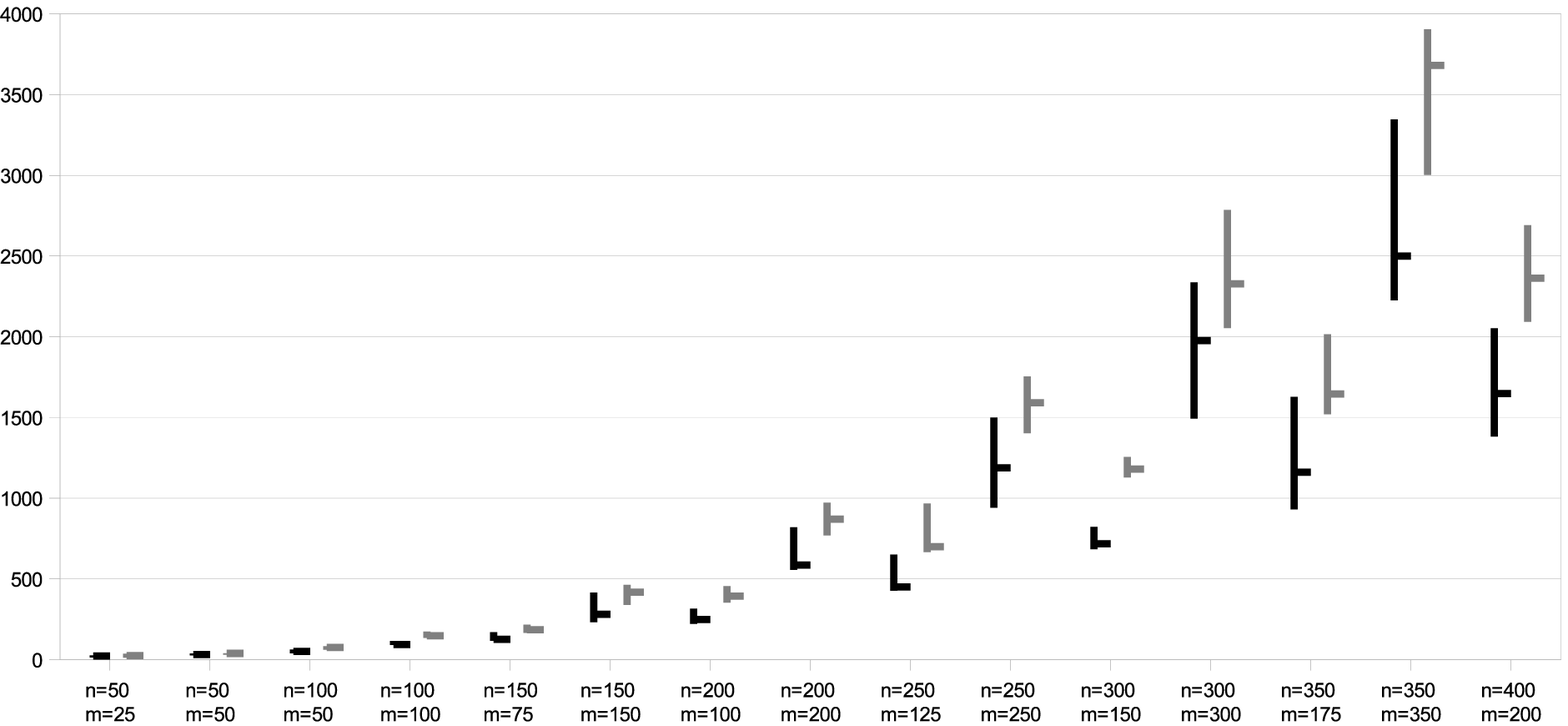}
    \captionof{figure}{Minimal, median and maximal solving time}
    \label{Figure:ComparisonSDproblemMedMinMax}
  \end{center}
\end{figure}
These results
justify that we will mainly concentrate on the reduced QCQP (\ref{SOCP:Luksan:Alg:QCQPTeilproblem2QsocpKapitel}) in the implementation as it is the only QCQP for which the solving time is competitive to that of the corresponding QP. 

\section{Numerical results}
\label{Paper:NumericalResults}
In the following section we compare the numerical results of our second order bundle algorithm with MPBNGC by \citet{MPBNGC} and SolvOpt by \citet{SolvOptPublication} for some examples of the Hock-Schittkowski collection by \citet{Schittkowski,Schittkowski2}, for custom examples that arise in the context of finding exclusion boxes for a quadratic CSP in GloptLab by \citet{Domes}, and for higher dimensional piecewise quadratic examples.
\subsection{Introduction}
There are three
implementations of Algorithm \ref{BundleSQPmitQCQP:Alg:GesamtAlgMitQCQP}
available:
A pure MATLAB version (for easy understanding, modifying and testing new ideas concerning the algorithm);
a MATLAB version in which the main parts of the algorithm are split into several subroutines, where every subroutine can either be called as pure MATLAB code or via a C \texttt{mex}-file (this is useful for partially speeding up the algorithm, but still keeping it simple enough for modifying and testing many examples of the modified code); and
a pure C version (for performance), which is used throughout all the tests.
The C \texttt{mex}-files and the C version require a BLAS/LAPACK implementation (e.g., ATLAS by \citet{ATLAS}, GotoBLAS by \citet{GotoBLAS}, or the Netlib BLAS reference implementation by \citet{NetlibBLAS}). In the unconstrained case, all three versions produce the same results as the original FORTRAN bundle-Newton method by \citet{Luksan}.

Although there exist some test collections for nonsmooth unconstrained optimization (e.g., \citet{Luksan3}) and nonsmooth linearly constrained optimization (e.g., \citet{Luksan2}; also cf.~\citet{KarmitsaComparison} for an extensive comparison of numerical results),
we do not know a standardized, prominent
test collection for nonsmooth constrained optimization. Therefore,
a
common way for testing nonsmooth constrained solvers is to take a test collection for smooth constrained optimization (e.g., the Hock-Schittkowski collection from \citet{Schittkowski,Schittkowski2}) and to treat the smooth constraints as one nonsmooth constraint (by using a $\max$-function). 

We will make tests for
\begin{itemize}
\item Algorithm \ref{BundleSQPmitQCQP:Alg:GesamtAlgMitQCQP} (with optimality tolerance $\varepsilon:=10^{-5}$), where we refer to the linearly constrained version as ``BNLC'', to the version with the QCQP (\ref{BundleSQPmitQCQP:Alg:QPTeilproblemReduced}) as
``Full Alg(orithm)'',
and to the version with the reduced QCQP (\ref{SOCP:Luksan:Alg:QCQPTeilproblem2QsocpKapitel}) as
``Red(uced) Alg(orithm)''
\item MPBNGC by \citet{MPBNGC} (with the standard termination criteria; although MPBNGC supports the handling of multiple nonsmooth constraints, we do not use this feature, since we are interested here, how well the different solvers handle the nonsmoothness of a constraint, i.e.~without exploiting the knowledge of the structure of a $\max$-function; since MPBNGC turned out to be very fast with respect to pure solving time for the low dimensional examples in the case of successful termination with a stationary point, the number of iterations and function evaluations was chosen in a way that in the other case the solving times of the different algorithms have approximately at least the same magnitude)
\item SolvOpt by \citet{SolvOptPublication} (with the standard termination criteria, which are described in \citet{SolvOpt})
\end{itemize}
(we choose MPBNGC and SolvOpt for our comparisons, since both are written in a compiled programming language, both are publicly available, and both support nonconvex constraints), where we will modify the termination criteria slightly only in
Subsection
\ref{section:HigherDimensionalPiecewiseQuadraticExamples}, on the following examples (the corresponding result tables can be found in
\ifthenelse{\boolean{AppendixOutsourcing}}
{\citet{HannesPaperAOutsourcing}):}
{Appendix \ref{section:ResultTables}):}
\begin{itemize}
\item
Optimization problem (\ref{BundleSQP:OptProblem}) with
$f(x):=\big(x_1+\tfrac{1}{2}\big)^2+\big(x_2+\tfrac{3}{2}\big)^2$
and
$F(x):=\max\hat{F}_{1:2}(x)$ (denoted by E1)
resp.~$F(x):=\max(-\hat{F}_{1:2}(x),\hat{F}_3(x))$ (denoted by E2),
where
$\hat{F}_1(x):= x_1^2+x_2^2-1$,
$\hat{F}_2(x):=(x_1-1)^2+(x_2+1)^2-1$, and
$\hat{F}_3(x):=(x_1-1)^2-x_2-1$,
the example from \citet[\GenaueAngabeFive]{HannesPaperB} (denoted by E3) and the Hock-Schittkowski collection (in the above sense; no problems which contain nonlinear equality constraints; linear constraints are inserted into the search direction problem in Algorithm \ref{BundleSQPmitQCQP:Alg:GesamtAlgMitQCQP}; feasible starting point). This yields 58 test problems 
which we will discuss in
Subsection
\ref{section:HStests}.
\item Optimization problems as described in \citet[\GenaueAngabeSix]{HannesPaperC} (for finding exclusion boxes for CSPs;
where the nonlinear part of these optimization problems is given by the certificate from \citet[\GenaueAngabeSeven]{HannesPaperC},
which we will discuss in
Subsection
\ref{section:ExclusionBoxTests}.
\item Higher dimensional piecewise quadratic examples with up to 100 variables
which we will discuss in
Subsection
\ref{section:HigherDimensionalPiecewiseQuadraticExamples}.
\end{itemize}
All test examples will be sorted with respect to the problem dimension (beginning with the smallest). Furthermore, we use analytic derivative information for all occurring functions (Note: Implementing analytic derivative information for the certificate from \citet[\GenaueAngabeEight]{HannesPaperC} effectively, is a nontrivial task) and we perform all tests on the same machine as in
Subsection
\ref{subsection:Comparisons}.

We introduce the following notation for the record of the solution process of an algorithm (which is used in this section as well as in
\ifthenelse{\boolean{AppendixOutsourcing}}
{\citet{HannesPaperAOutsourcing}).}
{Appendix \ref{section:ResultTables}).}
\begin{notation}
We define
\begin{align*}
\textnormal{N}
&:=
\textnormal{``Dimension of the optimization problem''}
\\
\textnormal{Nit}
&:=
\textnormal{``Number of performed iterations''}
\komma
\end{align*}
we denote the final number of evaluations of function dependent data by
\begin{align*}
\textnormal{Na}
&:=
\textnormal{``Number of calls to }(f,g,G,F,\hat{g},\hat{G})\textnormal{'' (Algorithm \ref{BundleSQPmitQCQP:Alg:GesamtAlgMitQCQP})}
\\
\NforMPBNGC
&:=
\textnormal{``Number of calls to }(f,g,F,\hat{g})\textnormal{'' (MPBNGC)}
\\
\NforSolvOpt
&:=
\textnormal{``Number of calls to }(f,F)\textnormal{'' (SolvOpt)}
\\
\textnormal{Ng}
&:=
\textnormal{``Number of calls to }g\textnormal{'' (SolvOpt)}
\\
\textnormal{N}\hat{\textnormal{g}}
&:=
\textnormal{``Number of calls to }\hat{g}\textnormal{'' (SolvOpt)}
\komma
\end{align*}
we denote the duration of the solution process by
\begin{align*}
\TotalTime
&:=
\textnormal{``Time in milliseconds''}
\\
\PartTime
&:=
\textnormal{``Time in milliseconds (without (QC)QP)'' (only relevant for Algorithm \ref{BundleSQPmitQCQP:Alg:GesamtAlgMitQCQP})}
\end{align*}
and we denote the additional algorithmic information by
\begin{align*}
\textnormal{R}
&:=
\textnormal{``Remark'' (e.g., if }t_0^k\textnormal{ is modified in Algorithm \ref{BundleSQPmitQCQP:Alg:GesamtAlgMitQCQP},}
\\
&\hspace{71pt}
\textnormal{additional SolvOpt termination information,}
\\
&\hspace{71pt}
\textnormal{supplementary problem dependent facts,\dots)}
\\
\textnormal{nt}
&:=
\textnormal{``No termination'' (within the given number of Nit,\dots)}
\\
\textnormal{wm}
&:=
\textnormal{``Wrong minimum''}
\punkt
\end{align*}
\end{notation}
\begin{remark}
In particular the percentage of the time spent in the (QC)QP in Algorithm \ref{BundleSQPmitQCQP:Alg:GesamtAlgMitQCQP} is given by
\begin{equation}
p_1
:=
\tfrac
{
\TotalTime(\textnormal{Algorithm \ref{BundleSQPmitQCQP:Alg:GesamtAlgMitQCQP}})
-
\PartTime(\textnormal{Algorithm \ref{BundleSQPmitQCQP:Alg:GesamtAlgMitQCQP}})
}{
\TotalTime(\textnormal{Algorithm \ref{BundleSQPmitQCQP:Alg:GesamtAlgMitQCQP}})
}
\punkt
\label{NumericalResults:Remark:p1}
\end{equation}
\end{remark}
For comparing the cost of evaluating function dependent data (like, e.g., function values, subgradients,\dots) in a
preferably
fair way (especially for solvers that use different function dependent data), we
will
make use of the following realistic ``credit point system'' that an optimal implementation of algorithmic differentiation in backward mode suggests (cf.~\citet{Griewank} and
\citet{schichl2004global,schichlcoconut,SchichlHabilitation}).
\begin{definition}
Let $f_A$, $g_A$ and $G_A$ resp.~$F_A$, $\hat{g}_A$ and $\hat{G}_A$ be the number of function values, subgradients and (substitutes of) Hessians of the objective function resp.~the constraint that an algorithm $A$ used for solving a nonsmooth optimization problem which may have linear constraints and at most one single nonsmooth nonlinear constraint. Then we define the cost of these evaluations by
\begin{equation}
c(A):=f_A+3g_A+3N\cdot G_A+\textnormal{nlc}\cdot(F_A+3\hat{g}_A+3N\cdot\hat{G}_A)
\label{NumericalResults:Definition:cost}
\komma
\end{equation}
where
$\textnormal{nlc}=1$ if the optimization problem has a nonsmooth nonlinear constraint, and $\textnormal{nlc}=0$ otherwise.
\end{definition}
Since Algorithm \ref{BundleSQPmitQCQP:Alg:GesamtAlgMitQCQP} evaluates $f$, $g$, $G$ and $F$, $\hat{g}$, $\hat{G}$ at every call that computes function dependent data, we obtain
\begin{equation*}
c(\textnormal{Algorithm \ref{BundleSQPmitQCQP:Alg:GesamtAlgMitQCQP}})
=
(1+\textnormal{nlc})\cdot\textnormal{Na}\cdot(1+3+3N)
\punkt
\end{equation*}
Since MPBNGC evaluates $f$, $g$ and $F$, $\hat{g}$ at every call that computes function dependent data (cf.~\citet{MPBNGC}), the only difference to Algorithm \ref{BundleSQPmitQCQP:Alg:GesamtAlgMitQCQP} with respect to $c$ from (\ref{NumericalResults:Definition:cost}) is that MPBNGC uses no information of Hessians and hence we obtain
\begin{equation*}
c(\textnormal{MPBNGC})
=
(1+\textnormal{nlc})\cdot\NforMPBNGC\cdot(1+3)
\punkt
\end{equation*}
Since SolvOpt evaluates $f$ and $F$ at every call that computes function dependent data and only sometimes $g$ or $\hat{g}$ (cf.~\citet{SolvOpt}), we obtain
\begin{equation*}
c(\textnormal{SolvOpt})
=
(1+\textnormal{nlc})\cdot\NforSolvOpt+3(\textnormal{Ng}+\textnormal{nlc}\cdot\textnormal{N}\hat{g})
\punkt
\end{equation*}

We will visualize the performance of two algorithms $A$ and $B$ for $s\in\lbrace c,\textnormal{Nit}\rbrace$ in
Subsection
\ref{section:HStests} and
Subsection
\ref{section:ExclusionBoxTests} by the following record-plot: In this plot the abscissa is labeled by the name of the test example and the value of the ordinate is given by
$
\rp(s):=s(B)-s(A)
$
(i.e.~if $\rp(s)>0$, then $\rp(s)$ tells us how much better algorithm $A$ is than algorithm $B$ with respect to $s$ for the considered example in absolute numbers; if $\rp(s)<0$, then $\rp(s)$ quantifies the advantage of algorithm $B$ in comparison to algorithm $A$; if $\rp(s)=0$, then both algorithms are equally good with respect to $s$). The scaling of the plots is chosen in a way that plots that contain the same test examples are comparable (although the plots may have been generated by results from different algorithms).
\begin{remark}\label{remark:socpBadMOSEKrobust}
All results for Algorithm \ref{BundleSQPmitQCQP:Alg:GesamtAlgMitQCQP} that are given in the tables of
\ifthenelse{\boolean{AppendixOutsourcing}}
{\citet{HannesPaperAOutsourcing})}
{Appendix \ref{section:ResultTables})}
were obtained by using MOSEK by \citet{AndersenPaper} for determining the search direction, where we used the MOSEK QCQP-solver which turned out to be much faster than the MOSEK SOCP-solver again (as we already noticed in Remark \ref{Remark:ComparisonsOfMOSEK:QCQPandSOCP}). We emphasize that in our tests there occurred no search direction problem which MOSEK was not able to solve.

The results for computing the search direction in Algorithm \ref{BundleSQPmitQCQP:Alg:GesamtAlgMitQCQP} with IPOPT by \citet{WaechterPaper} are practically the same with respect to Nit and Na. Furthermore, IPOPT was as robust and reliable as MOSEK. Nevertheless, IPOPT was slower than MOSEK with respect to the solving time which we expected as IPOPT is designed for general non-linear optimization problems,
while MOSEK is specialized in particular for QCQPs.

When using \texttt{socp} by \citet{BoydSOCP1995} for the computation of the search direction in Algorithm \ref{BundleSQPmitQCQP:Alg:GesamtAlgMitQCQP}, the results are also practically the same with respect to Nit and Na --- as long as \texttt{socp} did not fail to solve the search direction problem:
The most successful effort of stabilizing \texttt{socp} was achieved by the following idea from SEDUMI by \citet{Sturm,Polik}: We added an additional termination criterion to \texttt{socp} as it is used in SEDUMI, if SEDUMI cannot achieve the desired accuracy for the duality gap (the additional termination criterion is referred to as \texttt{pars.bigeps} in SEDUMI): If the current duality gap is smaller than $\texttt{bigeps}:=10^{-2}$ and differs at most by $10^{-5}$ from the duality gap of the last iteration, then we accept the current point as a solution.
In our empirical experiments \texttt{socp} tended to be more reliable, when we chose
certatin SOCP-dependent parameters according to \citet[\GenaueAngabeNEWtwo]{HannesDissertation}.
We were not able to make \texttt{socp} more robust by improving the strict feasibility of the starting point by solving various
linear programs
that are obtained from the primal SOCP
and the dual SOCP
by exploiting the fact that $\lvert x\rvert_{_2}\leq\lvert x\rvert_{_1}$ for all $x\in\mathbb{R}^n$ (\texttt{lp\_solve} by \citet{Berkelaar}, which is based on the revised simplex method and which we used for computing a solution of these linear programs, solved all of them easily).

At least when we used the variant of \texttt{socp} which was best for our purposes (i.e.~\texttt{socp} with a \texttt{bigeps}-termination criterion)
in Algorithm \ref{BundleSQPmitQCQP:Alg:GesamtAlgMitQCQP}, then we were able to solve all examples that we took from the Hock-Schittkowski collection, while we were not able to achieve this for the other variants of \texttt{socp}. Furthermore, many examples of the nonlinearly constrained optimization problem 
from \citet[\GenaueAngabeTen]{HannesPaperC}
were not solvable by Algorithm \ref{BundleSQPmitQCQP:Alg:GesamtAlgMitQCQP} when using \texttt{socp} for the computation of the search direction (even when we used the best variant of \texttt{socp}).
\end{remark}

\subsection{Hock-Schittkowski Test-set}
\label{section:HStests}
From Table 4
in
\ifthenelse{\boolean{AppendixOutsourcing}}
{\citet{HannesPaperAOutsourcing},}
{Appendix \ref{section:ResultTables},}
in which the results for the Hock-Schittkowski collection can be found and which is the basis for all plots in this subsection, we draw the following conclusions:

To compare the solving time $t_1$ for the reduced algorithm (with MOSEK as (QC)QP-solver) and MPBNGC, we consider
\begin{center}
\begin{normalsize}
\begin{tabular}{l|rrrr|}
                                  & \multicolumn{1}{c}{\TotalTime(Red Alg)}
                                  & \multicolumn{1}{c}{\PartTime(Red Alg)}
                                  & \multicolumn{1}{c}{$p_1$}
                                  & \multicolumn{1}{c|}{\TotalTime(MPBNGC)}
\\
\hline
HS     & 1198 & 961 & 0.80 & 1386\\
HS (*) &  902 & 751 & 0.83 &  154\\
\hline
\end{tabular}
\end{normalsize}
\end{center}
where we make use of
(\ref{NumericalResults:Remark:p1}) and in (*) we consider only those examples for which MPBNGC satisfied one of its termination criteria
(cf.~Subsubsection \ref{section:ExclusionBoxTests:NonLinearlyConstrainedCase}).
Hence, for those examples of the Hock-Schittkowski collection for which MPBNGC was able to terminate successfully, MPBNGC is faster than the reduced algorithm. Furthermore, we notice that the reduced algorithm spent at least 80\% of its time in the QCQP-solver, which is mostly overhead time in particular for the examples with lower dimension (which most examples are) as MOSEK has to, e.g., set up sparse matrix structures.

The reduced algorithm needs approximately 65\% of the solving time $t_1$ of the full algorithm. Nevertheless, SolvOpt only needs approximately 23\% resp.~36\% of the solving time $t_1$ of the full algorithm resp.~the reduced algorithm. Not surprisingly, the full algorithm spent 80\% of the time for solving the QCQPs (like the reduced algorithm did). Since SolvOpt terminated for the higher dimensional examples (i.e.~the 15-dimensional examples 284, 285 and 384) with points that are not stationary, while both the full and the reduced algorithm were able to solve them, and since the reduced algorithm needs significantly less pure solving time than the full algorithm for these examples
\begin{center}
\begin{normalsize}
\begin{tabular}{r|rrr|}
\multicolumn{1}{c|}{ex}
& \multicolumn{1}{c}{\TotalTime$($Full Alg$)$} & \multicolumn{1}{c}{\TotalTime$($Red Alg$)$} & \multicolumn{1}{c|}{$p_2$}\\
\hline
284 &   92 &  46 & 0.50\\
285 &  796 & 140 & 0.18\\
384 &  589 & 125 & 0.21\\
\hline
\end{tabular}
\end{normalsize}
\end{center}
where
$
p_2:=\tfrac{\TotalTime(\textnormal{Red Alg})}{\TotalTime(\textnormal{Full Alg})}
$,
we may expect that for more difficult examples the performance of the reduced algorithm increases with respect to $t_1$ (cf.~Subsubsection \ref{subsection:PureSolvingTime} and Subsection \ref{section:HigherDimensionalPiecewiseQuadraticExamples}).

Therefore, we will concentrate our comparison of Algorithm \ref{BundleSQPmitQCQP:Alg:GesamtAlgMitQCQP} (full and reduced version), MPBNGC and SolvOpt on the
qualitative aspects
of
the cost $c$ of the evaluations (solid line)
and
the number of iterations Nit (dashed line; this comparison is only meaningful for the comparison between the full algorithm and the reduced algorithm),
where we use the two different line types for a better distinction of the comparisons in Figure \ref{Figure:OverviewTableHS}, in this
subsection, where before making detailed comparisons of our 58 examples, we give a short overview of them as a reason of clarity of the presentation: This yields the following summary table consisting of the number of examples for which the reduced algorithm is better than the full algorithm, MPBNGC resp.~SolvOpt (and vice versa)
\begin{center}
\begin{scriptsize}
\begin{tabular}{l|rrrr|r|rrr|}
                                 & 
\multicolumn{1}{c}{no termi-}     & 
\multicolumn{1}{c}{significantly} & \multicolumn{1}{c}{better} & \multicolumn{1}{c}{a bit} &
\multicolumn{1}{|c|}{nearly} &
\multicolumn{1}{c}{a bit} & \multicolumn{1}{c}{better} & \multicolumn{1}{c|}{significantly}\\
                                 & 
\multicolumn{1}{c}{nation}     & 
\multicolumn{1}{c}{better} & \multicolumn{1}{c}{} & \multicolumn{1}{c}{better} &
\multicolumn{1}{|c|}{equal} &
\multicolumn{1}{c}{better} & \multicolumn{1}{c}{} & \multicolumn{1}{c|}{better}
\\
\hline
(Color code: Light grey)
& \multicolumn{4}{c|}{Full Alg} & & \multicolumn{3}{|c|}{Red Alg}\\
Nit             & 0 &	2	&	1	&	 3 & 51 &	 0 & 1 & 0\\
$c$             & 0 &	4	&	1	&	 7 & 40	&  1 & 1 & 4\\
\hline
(Color code: Grey)
& \multicolumn{4}{c|}{MPBNGC}  & & \multicolumn{3}{|c|}{Red Alg}\\
$c$             & 5 &	2	&	3	&	10 & 16	&  8 & 7 & 7\\
\hline
(Color code: Black)
& \multicolumn{4}{c|}{SolvOpt} & & \multicolumn{3}{|c|}{Red Alg}\\
$c$             & 3 &	3	&	4	&	 3 &	1	& 31 & 9 & 4\\
\hline
\end{tabular}
\end{scriptsize}
\end{center}
that is visualized in Figure \ref{Figure:OverviewTableHS}
\begin{center}
\captionsetup{type=figure}
\includegraphics[width=9cm]{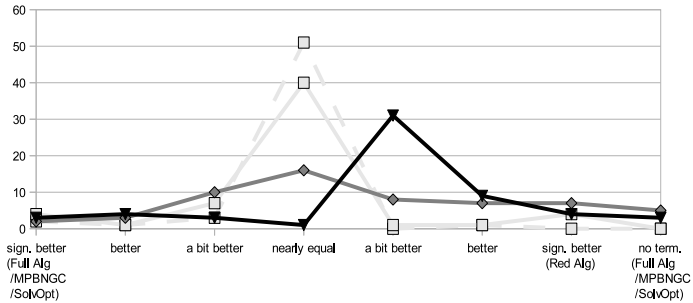}
\captionof{figure}{Hock-Schittkowski collection (summary)}
\label{Figure:OverviewTableHS}
\end{center}
and that let us draw the following conclusions:
The performances of the full algorithm and the reduced algorithm are quite similar.
The reduced algorithm is superior to MPBNGC in one third of the examples, for a further third of the examples one of these two solvers has only small advantages over the other, the performance differences between the two algorithms considered can be completely neglected for one quarter of the examples, and for the remaining ten percent of the examples MPBNGC beats the reduced algorithm clearly.
The reduced algorithm is superior to SolvOpt in about one quarter of the examples,
for sixty percent of the examples one of these two solvers has only small advantages over the other (in most cases the reduced algorithm is the slightly more successful one), and in the remaining twelve percent of the examples SolvOpt beats the reduced algorithm clearly.

Furthermore, only the full algorithm and the reduced algorithm solved all examples successfully.

\paragraph{Reduced algorithm vs. Full algorithm}
First of all, in the full algorithm $t_0^k$ is only modified in 11 examples (34, 43, 66, 83, 100, 113, 227, 230, 264, 285, 384), while in the reduced algorithm this happens in 14 examples (the additional examples are 284, 330, 341). In all these examples $t_0^k$ is only modified a few times and a modification only occurs at very early iterations of the optimization process (cf.~\citet[\GenaueAngabeEleven]{HannesPaperB}).

From Figure 14 
and Figure 15 in \citet{HannesPaperAOutsourcing}
we conclude that the full and the reduced algorithm produce in most of the 58 examples approximately the same results --- exceptions from this observation are in view of iterations the following
7 examples: The reduced algorithm is
better in 1 example
in comparison with the full algorithm, while the full algorithm is
significantly better in 2 examples,
better in 1 example
and
a bit better in 3 examples
in comparison with the reduced algorithm.

In view of costs the exceptions are given by the following
18 examples: The reduced algorithm is
significantly better in 4 examples,
better in 1 example (33)
a bit better in 1 example
in comparison with the full algorithm, while the full algorithm is
significantly better in 4 examples,
better in 1 example
and
a bit better in 7 examples
in comparison with the reduced algorithm.

\paragraph{Reduced algorithm vs. MPBNGC}
MPBNGC does not satisfy any of its termination criteria for five examples (15, 20, 83, 285 and 384) within the given number of iterations and function evaluations.
For the other 53 examples from Figure 16 in \citet{HannesPaperAOutsourcing}
we emphasize the following ones: The reduced algorithm is
significantly better in 7 examples,
better in 7 examples
and
a bit better in 8 examples
in comparison with MPBNGC, while MPBNGC is
significantly better in 2 examples,
better in 3 examples
and
a bit better in 10 examples
in comparison with the reduced algorithm. In the remaining
16 examples the cost of the reduced algorithm and MPBNGC is practically the same.

\paragraph{Reduced algorithm vs. SolvOpt}
SolvOpt terminates for the three 15-dimensional examples 284, 285 and 384 with points that are not stationary.
For the other 55 examples from Figure 17 in \citet{HannesPaperAOutsourcing}
we emphasize the following ones: The reduced algorithm is
significantly better in 4 examples
and
better in 9 examples
in comparison with SolvOpt, while SolvOpt is
significantly better in 3 examples,
better in 4 examples
and
a bit better in 3 examples
in comparison with the reduced algorithm. Except for example 233 in which the cost of the reduced algorithm and SolvOpt are practically the same, in all
31 remaining examples the reduced algorithm is a bit better than SolvOpt.

\subsection{Exclusion boxes}
\label{section:ExclusionBoxTests}
\ifthenelse{\boolean{ResultTablesBlack}}
{
\subsubsection{Basics}
}
{
\subsubsection{Notation}
}
\label{subsection:ExclusionBoxTestsNotation}
We consider the
quadratic
CSP
\begin{equation}
\begin{split}
F(x)&\in\boldsymbol{F}\\
  x &\in\boldsymbol{x}
\label{Abstract:coi:CSP}
\end{split}
\end{equation}
and we assume that a solver, which is able to solve a CSP,
takes the box $\boldsymbol{u}:=[\underline{u},\overline{u}]\subseteq\boldsymbol{x}$ into consideration during the solution process.
\citet[\GenaueAngabeCertificate]{HannesPaperC}
constructed a certificate of infeasibility $f$, which is a nondifferentiable and nonconvex function in general, with the following property: If there exists a vector $y$ with
\begin{equation}
f(y,\underline{u},\overline{u})<0\komma
\label{Abstract:CertificateNegative}
\end{equation}
then the CSP (\ref{Abstract:coi:CSP}) has no feasible point in $\boldsymbol{u}$ and consequently this box can be excluded for the rest of the solution process. Therefore,
a box $\boldsymbol{u}$ for which (\ref{Abstract:CertificateNegative}) holds is called an \textbf{exclusion box}.

The obvious way for finding an exclusion box for the CSP (\ref{Abstract:coi:CSP}) is to minimize $f$
\begin{equation*}
\begin{split}
\min_{y}{f(y,\underline{u},\overline{u})}
\end{split}
\end{equation*}
and stop the minimization if a negative function value occurs. We
will
give results for
this
linearly constrained optimization problem
with a fixed box (i.e.~without optimizing $u$ and $v$) for dimensions between $4$ and $11$ in
Subsubsection
\ref{section:ExclusionBoxTests:LinearlyConstrainedCaseFixedInterval}.

To find at least an exclusion box $\boldsymbol{v}:=[\underline{v},\overline{v}]\subseteq\boldsymbol{u}$ with $\underline{v}+r\leq\overline{v}$, where $r\in(\zeroVector{n},\overline{u}-\underline{u})$ is fixed, we can try to solve 
\begin{equation*}
\begin{split}
&\min_{y,\underline{v},\overline{v}}{f(y,\underline{v},\overline{v})}\\
&\textnormal{ s.t. }[\underline{v}+r,\overline{v}]\subseteq\boldsymbol{u}
\komma
\end{split}
\end{equation*}
where
the
results for
this
linearly constrained optimization problem
with a variable box (i.e.~with optimizing $u$ and $v$) for dimensions between $8$ and $21$
are discussed
in
Subsubsection
\ref{section:ExclusionBoxTests:LinearlyConstrainedCaseVariableInterval}.

Moreover, we can enlarge an exclusion box $\boldsymbol{v}$ by solving
\begin{equation*}
\begin{split}
&\max_{y,\underline{v},\overline{v}}{\mu(\underline{v},\overline{v})}\\
&\textnormal{ s.t. }f(y,\underline{v},\overline{v})\leq\delta\komma~[\underline{v},\overline{v}]\subseteq\boldsymbol{u}\komma
\end{split}
\end{equation*}
where $\delta<0$ is given and
$
\mu(\underline{v},\overline{v})
:=
\lvert
\left(
\begin{smallmatrix}
\underline{v}-\underline{x}\\
\overline{v}-\overline{x}
\end{smallmatrix}
\right)
\rvert_{_1}
$
measures the magnitude of the box $\boldsymbol{v}$,
and we regard an exclusion box as sufficiently large, if the objective function satisfies
$\mu(\underline{v},\overline{v})\leq10^{-6}$.
The
discussion of the results of
this
nonlinearly constrained optimization problem
for dimension $8$
can be found in
in
Subsubsection
\ref{section:ExclusionBoxTests:NonLinearlyConstrainedCase}.

The underlying data for these nonsmooth optimization problems was extracted from real CSPs that occur in GloptLab by \citet{Domes}. Apart from $u$ and $v$, we will concentrate on the optimization of the variables $y$ and $z$ due to the large number of tested examples (cf.~Subsubsection \ref{subsection:PureSolvingTime}), and since the additional optimization of $R$ and $S$ did not have much impact on the quality of the results which was discovered in additional empirical observations, where a detailed analysis of these observations goes beyond the scope of this
paper.
Furthermore, we will make our tests for the two different choices $T=1$ and $T=\lvert y\rvert_{_2}$ of the function $T$, which occurs in the denominator of the certificate $f$ from \citet[\GenaueAngabeFifteen]{HannesPaperC}, where for the latter one $f$ is only defined outside of the zero set of $T$ which has measure zero --- although the convergence theory of many solvers (cf., e.g., \citet[\GenaueAngabeConvergenceTheory]{HannesPaperB})
requires that all occurring functions are defined on the whole
space.
\ifthenelse{\boolean{ResultTablesBlack}}
{
}
{

We represent the information on the feasibility of example $i$ that is obtained by GloptLab in the following way
\ifthenelse{\boolean{AppendixOutsourcing}}
{(cf.~\citet{HannesPaperAOutsourcing}):}
{(cf.~Appendix \ref{section:ResultTables}):}
\begin{center}
\begin{tabular}{l|l|}
\multicolumn{1}{c|}{ex} & \multicolumn{1}{c|}{Status of example $i$ as identified by GloptLab}\\
\hline
\textcolor{red}{$i$}   & infeasible\\
\textcolor{blue}{$i$}  & feasible\\ 
$i$                    & No statement on (in)feasibility can be made\\
\textcolor{green}{$i$} & A local minimizer is found\\
\textcolor{gold}{$i$}  & Maximum number of iterations is exceeded\\
\hline
\end{tabular}
\end{center}
We represent the status of a point concerning certificate $f$ from \citet[\GenaueAngabeFifteen]{HannesPaperC} after performing Nit iterations of a nonsmooth solver in the linearly constrained case 
(cf.~Subsubsection~\ref{section:ExclusionBoxTests:LinearlyConstrainedCaseFixedInterval} and Subsubsection \ref{section:ExclusionBoxTests:LinearlyConstrainedCaseVariableInterval}) by
\begin{center}
\begin{tabular}{l|l|}
\multicolumn{1}{c|}{Nit} & \multicolumn{1}{c|}{Status after Nit iterations (linearly constrained case)}\\
\hline
\textcolor{red}{k}    & A point with $f<0$ is found\\
                k     & A stationary point with $f\geq0$ is found\\
\textcolor{blue}{k}   & A feasible point is found (of course with $f\geq0$)\\
\textcolor{gold}{$k$} & Maximum number of iterations is exceeded (without (in)feasibility statement)\\
\hline
\end{tabular}
\end{center}
and in the nonlinearly constrained case
(cf.~Subsubsection~\ref{section:ExclusionBoxTests:NonLinearlyConstrainedCase})
by
\begin{center}
\begin{tabular}{l|l|}
\multicolumn{1}{c|}{Nit} & \multicolumn{1}{c|}{Status after Nit iterations (nonlinearly constrained case)}\\
\hline
\textcolor{red}{k}    & A point which yields a sufficiently large exclusion box is found\\
                k     & A stationary point is found where the exclusion box is not large enough\\
\textcolor{gold}{$k$} & Maximum number of iterations is exceeded (without finding a sufficiently\\
                      & large exclusion box or being able to make a stationarity statement)\\
\hline
\end{tabular}
\end{center}
where an exclusion box is sufficiently large in the nonlinearly constrained case, if the objective function satisfies
\begin{equation*}
b_+^1(y,z,R,S,u,v)\leq10^{-6}\punkt
\end{equation*}
}
\begin{remark}
Because SolvOpt cannot distinguish between linear and nonlinear constraints (cf.~\citet[p.~15]{SolvOpt}), the linear constraints of the linearly constrained optimization problems
from \citet[\GenaueAngabeSixteen]{HannesPaperC}
must be formulated as nonlinear constraints in SolvOpt. Nevertheless, we will not include the number of these evaluations in the computation of the cost $c$ from (\ref{NumericalResults:Definition:cost}) for the mentioned optimization problems in
Subsubsection
\ref{section:ExclusionBoxTests:LinearlyConstrainedCaseFixedInterval} and
Subsubsection
\ref{section:ExclusionBoxTests:LinearlyConstrainedCaseVariableInterval}, since these evaluations may be considered as easy in comparison to the evaluation of the certificate $f$ from
\citet[\GenaueAngabeSeventeen]{HannesPaperC}
which is the objective function in these optimization problems.
\end{remark}

\subsubsection{Overview of the results}
\label{subsection:PureSolvingTime}
We compare the total time $t_1$ of the solution process, where we used the reduced algorithm (with MOSEK as
the
(QC)QP-solver) in the constrained case: From
Tables 5--8
\ifthenelse{\boolean{AppendixOutsourcing}}
{(s.~\citet{HannesPaperAOutsourcing})}
{(s.~Appendix \ref{section:ResultTables})}
we obtain
\begin{center}
\begin{footnotesize}
\begin{tabular}{l|rrrrr|}
                                  & \multicolumn{1}{c}{\TotalTime(Red Alg)}
                                  & \multicolumn{1}{c}{\PartTime(Red Alg)}
                                  & \multicolumn{1}{c}{$p_1$}
                                  & \multicolumn{1}{c}{\TotalTime(MPBNGC)}
                                  & \multicolumn{1}{c|}{\TotalTime(SolvOpt)}
\\
\hline
                                     & \multicolumn{5}{c|}{$T=1$}\\
Linearly constrained (fixed box)     &  1477 &  215 & 0.85 &   231 &  2754 \\
Linearly constrained (variable box)  &   782 &   60 & 0.92 &    30 &  1546 \\
Nonlinearly constrained              & 25420 & 4885 & 0.81 & 21860 &  38761\\
Nonlinearly constrained (*)          & 19053 & 3723 & 0.80 &  2067 &  30312\\
\hline
                                     & \multicolumn{5}{c|}{$T=\lvert y\rvert_{_2}$}\\
Linearly constrained (fixed box)     &  1316 &  129 & 0.90 &    15 &  1508\\
Linearly constrained (variable box)  &   797 &   45 & 0.94 &    30 &  2263\\
Nonlinearly constrained              & 24055 & 4284 & 0.82 & 25383 & 16909\\
Nonlinearly constrained (*)          & 18038 & 3112 & 0.83 &  3719 & 12635\\
\hline
\end{tabular}
\end{footnotesize}
\end{center}
where we make use of
(\ref{NumericalResults:Remark:p1}) and in (*) we consider only those examples for which MPBNGC satisfied one of its termination criteria
(cf.~Subsubsection \ref{section:ExclusionBoxTests:NonLinearlyConstrainedCase}).

For the linearly constrained problems MPBNGC was the fastest of the tested algorithms, followed by BNLC and SolvOpt.
If we consider only those nonlinearly constrained examples for which MPBNGC was able to terminate successfully, MPBNGC was the fastest algorithm again. Considering the competitors, for the nonlinearly constrained problems with $T=1$ the reduced algorithm is 13.3 seconds resp.~11.3 seconds faster than SolvOpt, while for the nonlinearly constrained problems with $T=\lvert y\rvert_{_2}$ SolvOpt is 7.1 seconds resp.~5.4 seconds faster than the reduced algorithm.

Again (cf.~Subsection \ref{section:HStests}), taking a closer look at $p_1$ yields the observation that at least 85\% of the time is consumed by solving the QP (in the linearly constrained case) resp.~at least 80\% of the time is consumed by solving the QCQP (in the nonlinearly constrained case), which implies that the difference in the percentage between the QP and the QCQP is small in particular (an investigation of the behavior of the solving time $t_1$ for higher dimensional problems can be found in
Subsection
\ref{section:HigherDimensionalPiecewiseQuadraticExamples}).

Therefore, we will concentrate in
Subsubsection
\ref{section:ExclusionBoxTests:LinearlyConstrainedCaseFixedInterval},
Subsubsection
\ref{section:ExclusionBoxTests:LinearlyConstrainedCaseVariableInterval} and
Subsubsection
\ref{section:ExclusionBoxTests:NonLinearlyConstrainedCase} on the comparison of qualitative aspects between Algorithm \ref{BundleSQPmitQCQP:Alg:GesamtAlgMitQCQP}, MPBNGC and SolvOpt (like, e.g., the cost $c$ of the evaluations),
where before making these detailed comparisons, we give a short overview of them as a reason of clarity of the presentation: In both cases
$T=1$ (solid line) and
$T=\lvert y\rvert_{_2}$ (dashed line),
where we use the two different line types for a better distinction in the following, we tested
128 linearly constrained examples with a fixed box,
117 linearly constrained examples with a variable box and
201 nonlinearly constrained examples,
which yields the following two summary tables consisting of the number of examples for which Algorithm \ref{BundleSQPmitQCQP:Alg:GesamtAlgMitQCQP} (BNLC resp.~the reduced algorithm) is better than MPBNGC resp.~SolvOpt (and vice versa) with respect to the cost $c$ of the evaluations
\begin{center}
\begin{tiny}
\begin{tabular}{l|rrrrrrrr|}
(Color code: Light grey)
& \multicolumn{4}{c|}{MPBNGC} & & \multicolumn{3}{|c|}{BNLC/Red Alg}\\
                                 &
\multicolumn{1}{c}{no termi-}     & 
\multicolumn{1}{c}{significantly} & \multicolumn{1}{c}{better} & \multicolumn{1}{c}{a bit} &
\multicolumn{1}{|c|}{nearly} &
\multicolumn{1}{c}{a bit} & \multicolumn{1}{c}{better} & \multicolumn{1}{c|}{significantly}
\\
                                 &
\multicolumn{1}{c}{nation}     & 
\multicolumn{1}{c}{better} & \multicolumn{1}{c}{} & \multicolumn{1}{c}{better} &
\multicolumn{1}{|c|}{equal} &
\multicolumn{1}{c}{better} & \multicolumn{1}{c}{} & \multicolumn{1}{c|}{better}
\\
\hline
                                     &    &  &  &  & \multicolumn{1}{c}{$T=1$} & & &\\
Linearly constrained (fixed box)     &  0 & 2 & 5  & 12 & 106 &  2 & 0 & 1\\
Linearly constrained (variable box)  &  0 & 0 & 0  &  1 & 116 &  0 & 0 & 0\\
Nonlinearly constrained              & 32 & 6 & 28 & 89 &  31 & 10 & 2 & 3\\
\hline
                                     &  &  &  &  & \multicolumn{1}{c}{$T=\lvert y\rvert_{_2}$} & & &\\
Linearly constrained (fixed box)     &  0 & 2 & 5 &  30 &  91 &  0 &  0 & 0\\
Linearly constrained (variable box)  &  0 & 0 & 0 &   5 & 112 &  0 &  0 & 0\\
Nonlinearly constrained              & 43 & 4 & 28 & 59 &  30 & 15 & 14 & 8\\
\hline
\end{tabular}
\end{tiny}
\end{center}
\begin{center}
\begin{tiny}
\begin{tabular}{l|rrrrrrrr|}
(Color code: Black)
& \multicolumn{4}{c|}{SolvOpt} & & \multicolumn{3}{|c|}{BNLC/Red Alg}\\
                                 &
\multicolumn{1}{c}{no termi-}     & 
\multicolumn{1}{c}{significantly} & \multicolumn{1}{c}{better} & \multicolumn{1}{c}{a bit} &
\multicolumn{1}{|c|}{nearly} &
\multicolumn{1}{c}{a bit} & \multicolumn{1}{c}{better} & \multicolumn{1}{c|}{significantly}
\\
                                 &
\multicolumn{1}{c}{nation}     & 
\multicolumn{1}{c}{better} & \multicolumn{1}{c}{} & \multicolumn{1}{c}{better} &
\multicolumn{1}{|c|}{equal} &
\multicolumn{1}{c}{better} & \multicolumn{1}{c}{} & \multicolumn{1}{c|}{better}
\\
\hline
                                     &  &  &  &  & \multicolumn{1}{c}{$T=1$} & & &\\
Linearly constrained (fixed box)     & 0 & 1 &  3 &  0 & 61 & 25 & 13 & 25\\
Linearly constrained (variable box)  & 0 & 0 &  0 &  0 & 48 & 37 & 24 & 8\\
Nonlinearly constrained              & 0 & 0 & 14 & 20 & 21 & 76 & 20 & 50\\
\hline
                                     &  &  &  &  & \multicolumn{1}{c}{$T=\lvert y\rvert_{_2}$} & & &\\
Linearly constrained (fixed box)     & 0 & 1 &  2 &  1 & 32 & 34 & 49 &  9\\
Linearly constrained (variable box)  & 0 & 0 &  0 &  5 & 41 & 32 & 19 & 20\\
Nonlinearly constrained              & 0 & 2 & 24 & 26 & 31 & 61 & 45 & 12\\
\hline
\end{tabular}
\end{tiny}
\end{center}
that are visualized in
Figures \ref{Figure:OverviewTableLinearlyConstrainedFixedBox},
\ref{Figure:OverviewTableLinearlyConstrainedVariableBox}, and
\ref{Figure:OverviewTableNonlinearlyConstrained}
\\
\begin{center}
\captionsetup{type=figure}
\includegraphics[width=8cm]{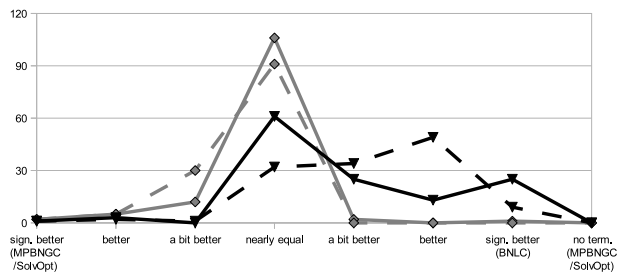}
\captionof{figure}{Linearly constrained --- fixed box (summary)}
\label{Figure:OverviewTableLinearlyConstrainedFixedBox}
\end{center}
\begin{center}
\captionsetup{type=figure}
\includegraphics[width=8cm]{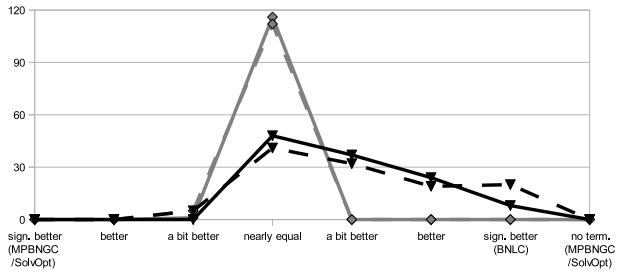}
\captionof{figure}{Linearly constrained --- variable box (summary)}
\label{Figure:OverviewTableLinearlyConstrainedVariableBox}
\end{center}
\begin{center}
\captionsetup{type=figure}
\includegraphics[width=8cm]{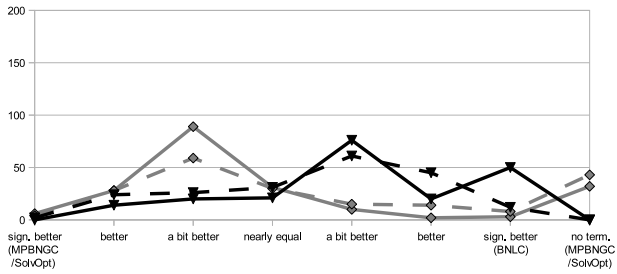}
\captionof{figure}{Nonlinearly constrained (summary)}
\label{Figure:OverviewTableNonlinearlyConstrained}
\end{center}
and that let us draw the following conclusions:

The performance differences between BNLC and MPBNGC can be neglected for the largest part of the linearly constrained examples (with small advantages for MPBNGC in about ten percent of these examples).
For the nonlinearly constrained examples the reduced algorithm is superior to MPBNGC in one quarter of the examples, for forty percent of the examples one of these two solvers has small advantages over the other (in most cases MPBNGC is the slightly more successful one), the performance differences between the two algorithms considered can be completely neglected for fifteen percent of the examples, and for further fifteen percent of the examples MPBNGC beats the reduced algorithm clearly.

For the linearly constrained examples BNLC is superior to SolvOpt in one third of the examples, for one quarter of the examples one of these two solvers has small advantages over the other (in nearly all cases BNLC is the slightly more successful one), the performance differences between the two algorithms considered can be completely neglected for forty percent of the examples, and in only one percent of the examples SolvOpt beats the reduced algorithm clearly.
For the nonlinearly constrained examples the reduced algorithm is superior to SolvOpt in one third of the examples, for 45 percent of the examples one of these two solvers has small advantages over the other (the reduced algorithm is often the slightly more successful one), the performance differences between the considered two algorithms can be completely neglected for ten percent of the examples, and in the remaining ten percent of the examples SolvOpt beats the reduced algorithm clearly.

In contrast to the linearly constrained case, in which all three solvers terminated successfully for all examples, only the reduced algorithm and SolvOpt were able to attain this goal in the nonlinearly constrained case, too.

\subsubsection{Linearly constrained case (fixed box)}
\label{section:ExclusionBoxTests:LinearlyConstrainedCaseFixedInterval}
We took 310 examples from real CSPs that occur in GloptLab. We observe that
for 79 examples
the starting point is feasible for the CSP and
for 103 examples
the evaluation of the certificate at the starting point identifies the box as infeasible
and hence there remain 128 test problems.

\paragraph{BNLC vs. MPBNGC}
In the case $T=1$ we conclude from Figure 18 in \citet{HannesPaperAOutsourcing}
that BNLC is
significantly better in 1 example
and
a bit better in 2 examples
in comparison with MPBNGC, while MPBNGC is
significantly better in 2 examples,
better in 5 examples
and
a bit better in 12 examples
in comparison with BNLC. In the 106 remaining examples the costs of BNLC and MPBNGC are practically the same.

In the case $T=\lvert y\rvert_{_2}$ it follows from Figure 19 in \citet{HannesPaperAOutsourcing}
that MPBNGC is
significantly better in 2 examples,
better in 5 examples
and
a bit better in 30 examples
in comparison with BNLC. In the 91 remaining examples the costs of BNLC and MPBNGC are practically the same.

\paragraph{BNLC vs. SolvOpt}
In the case $T=1$ we conclude from Figure 20 in \citet{HannesPaperAOutsourcing}
that BNLC is
significantly better in 25 examples,
better in 13 examples
and
a bit better in 25 examples
in comparison with SolvOpt, while SolvOpt is
significantly better in 1 example
and
better in 3 examples
in comparison with BNLC. In the 61 remaining examples the costs of BNLC and SolvOpt are practically the same.

In the case $T=\lvert y\rvert_{_2}$ it follows from Figure 21 in \citet{HannesPaperAOutsourcing}
that BNLC is
significantly better in 9 examples,
better in 49 examples
and
a bit better in 34 examples
in comparison with SolvOpt, while SolvOpt is
significantly better in 1 example,
better in 2 examples
and
a bit better in 1 example
in comparison with BNLC. In the 32 remaining examples the costs of BNLC and SolvOpt are practically the same.

\subsubsection{Linearly constrained case (variable box)}
\label{section:ExclusionBoxTests:LinearlyConstrainedCaseVariableInterval}
We observe that
for 80 examples
the starting point is feasible for the CSP and
for 113 examples
the evaluation of the certificate at the starting point identifies the boxes as infeasible
and hence there remain 117 test problems of the 310 original examples from GloptLab.

\paragraph{BNLC vs. MPBNGC}
In the case $T=1$ we conclude from Figure 22 in \citet{HannesPaperAOutsourcing}
that MPBNGC is
a bit better in 1 example
in comparison with BNLC. In the 116 remaining examples the costs of BNLC and MPBNGC are practically the same.

In the case $T=\lvert y\rvert_{_2}$ it follows from Figure 23 in \citet{HannesPaperAOutsourcing}
that MPBNGC is
a bit better in 5 examples
in comparison with BNLC. In the 112 remaining examples the costs of BNLC and MPBNGC are practically the same.

\paragraph{BNLC vs. SolvOpt}
In the case $T=1$ we conclude from Figure 24 in \citet{HannesPaperAOutsourcing}
that BNLC is
significantly better in 8 examples,
better in 24 examples
and
a bit better in 37 examples
in comparison with SolvOpt. In the 48 remaining examples the costs of BNLC and SolvOpt are practically the same.

In the case $T=\lvert y\rvert_{_2}$ it follows from Figure 25 in \citet{HannesPaperAOutsourcing}
that BNLC is
significantly better in 20 examples,
better in 19 examples
and
a bit better in 32 examples
in comparison with SolvOpt, while SolvOpt is
a bit better in 5 examples
in comparison with BNLC. In the 41 remaining examples the costs of BNLC and SolvOpt are practically the same.

\subsubsection{Nonlinearly constrained case}
\label{section:ExclusionBoxTests:NonLinearlyConstrainedCase}
Since we were not able to find a starting point, i.e.~an infeasible sub-box, for 109 examples,
we exclude them from the following tests for which there remain 201 examples of the 310 original examples from GloptLab.

\paragraph{Reduced algorithm vs. MPBNGC}
In the case $T=1$ MPBNGC does not satisfy any of its termination criteria for 32 examples
within the given number of iterations and function evaluations
(also cf.~Subsubsection \ref{subsection:ExclusionBoxTestsNotation}).
For the remaining
169 examples we conclude from Figure 26 in \citet{HannesPaperAOutsourcing}
that the reduced algorithm is
significantly better in 3 examples,
better in 2 examples
and
a bit better in 10 examples
in comparison with MPBNGC, while MPBNGC is
significantly better in 6 examples,
better in 28 examples
and
a bit better in 89 examples
in comparison with the reduced algorithm, and in
31
examples the costs of the reduced algorithm and MPBNGC are practically the same.

In the case $T=\lvert y\rvert_{_2}$ MPBNGC does not satisfy any of its termination criteria for 43 examples
within the given number of iterations and function evaluations. For the remaining
158
examples it follows from Figure 27 in \citet{HannesPaperAOutsourcing}
that the reduced algorithm is
significantly better in 8 examples,
better in 14 examples
and
a bit better in 15 examples
in comparison with MPBNGC, while MPBNGC is
significantly better in 4 examples,
better in 28 examples
and
a bit better in 59 examples
in comparison with the reduced algorithm, and in
30
examples the costs of the reduced algorithm and MPBNGC are practically the same.

\paragraph{Reduced algorithm vs. SolvOpt}
In the case $T=1$ we conclude from Figure 28 in \citet{HannesPaperAOutsourcing}
that the reduced algorithm is
significantly better in 50 examples,
better in 20 examples
and
a bit better in 76 examples
in comparison with SolvOpt, while SolvOpt is
better in 14 examples
and
a bit better in 20 examples
in comparison with the reduced algorithm. In the 21 remaining examples the costs of the reduced algorithm and SolvOpt are practically the same.

In the case $T=\lvert y\rvert_{_2}$ it follows from Figure 29 in \citet{HannesPaperAOutsourcing}
that the reduced algorithm is
significantly better in 12 examples,
better in 45 examples
and
a bit better in 61 examples
in comparison with SolvOpt, while SolvOpt is
significantly better in 2 examples,
better in 24 examples
and
a bit better in 26 examples
in comparison with the reduced algorithm. In the 31 remaining examples the costs of the reduced algorithm and SolvOpt are practically the same.

\subsection{Higher dimensional piecewise quadratic examples}
\label{section:HigherDimensionalPiecewiseQuadraticExamples}
We want to give numerical results for the nonsmooth optimization problem (\ref{BundleSQP:OptProblem}) with
\begin{equation*}
f(x)
:=
\max_{i=1,\dots,m_1}{f_i(x)}
\komma\quad
F(x)
:=
\max_{j=1,\dots,m_2}{F_j(x)}
\komma
\end{equation*}
where
\begin{align*}
f_i(x)
&:=
\alpha_i+a_i^T(x-x_i)+\tfrac{1}{2}(x-x_i)^TA_i(x-x_i)
\\
F_j(x)
&:=
\beta_j+b_j^T(x-x_j)+\tfrac{1}{2}(x-x_j)^TB_j(x-x_j)
\end{align*}
and $\alpha_i,\beta_j\in\mathbb{R}$, $a_i,b_j\in\mathbb{R}^N$, $A_i,B_j\in\Sym{N}$, $x_i,x_j\in\mathbb{R}^N$.

The underlying data of the test examples was produced by a random number generator with the following restrictions concerning the data corresponding to $F$:
At least one $B_j$ is chosen as a positive definite matrix to guarantee that the feasible set is bounded,
and after
choosing $b_j$, $B_j$, $x_j$ as well as a starting point $x^0\in\mathbb{R}^N$, $\beta_j$ is chosen such that $x^0$ is strictly feasible.

We made comparison tests for the dimensions
$N\in\lbrace20,40,60,80,100\rbrace$ to investigate the behavior of
the reduced
algorithm ($\square$),
MPBNGC ($\diamondsuit$) and
SolvOpt ($\nabla$),
where we use the colors to distinguish the results of the different solvers,
with respect to
the
solving time $t_1$
and
successful termination,
and we focus
on the larger values of $N$ (due to the magnitude of $N$, we did not test the full version of Algorithm \ref{BundleSQPmitQCQP:Alg:GesamtAlgMitQCQP}). Moreover, we chose $m_1:=\tfrac{N}{10}$ and $m_2\in\lbrace\tfrac{N}{2},N\rbrace$, so that the emphasis of the examples lies on the handling of the constraint.

Furthermore, due to the magnitude of the test examples, we weakened the optimality tolerance of the reduced algorithm to $\varepsilon:=10^{-3}$. Since the reduced algorithm terminated for all examples of this class of test functions with satisfying its termination criterion (which guarantees the stationarity of the computed point due to
\citet[\GenaueAngabeHD]{HannesPaperB}), we denote the minimizer (of the corresponding example) that was computed by the reduced algorithm by $\hat{x}$.

Before the actual tests, we performed a few runs of the whole test set, where we started with very weak termination criteria for MPBNGC and SolvOpt and then sharpened them, with the goal to make the results between the different solvers comparable in the following way: If the computed minimizer is close to $\hat{x}$, then approximately the same $F_j$ should be active. Based on these empirical observations, we made the final choices for the termination criteria of MPBNGC and SolvOpt, where we were quite successful to achieve this goal for MPBNGC, while we were not able to achieve it for SolvOpt in many cases (although putting a lot of effort into it).

For every pair $(N,m_2)$ we tested $20$ different examples for two levels of difficulty that is classified by the average number of $j\in\lbrace1,\dots,m_2\rbrace$ with
$\lvert F_j(\hat{x})-F(\hat{x})\rvert
\leq
10^{-3}$,
which yields the following overview of our overall 400 different examples
\begin{center}
\begin{normalsize}
\begin{tabular}{l|l|rrrrr|}
\multicolumn{1}{c|}{Level} & \multicolumn{1}{c|}{$m_2$} & \multicolumn{5}{c|}{$N$}\\
                           &                            & 20 & 40 & 60 &80 & 100\\
\hline
Easy                       & $\tfrac{N}{2}$             &  4 &  4 &  6 &  6 &   7\\
                           & $N$                        &  5 &  6 &  8 &  9 &  10\\
\hline
Difficult                  & $\tfrac{N}{2}$             &  4 &  8 & 12 & 15 &  19\\
                           & $N$                        &  7 & 14 & 19 & 26 &  31\\
\hline
\end{tabular}
\end{normalsize}
\end{center}
i.e.~for given $N$ and $m_2$ we regard an example as more difficult,
the more impact the constraint has at $\hat{x}$ (in the case of the
successful termination of one of the solvers, there was always at
least one $F_j$ active). Moreover, for a given level of difficulty,
$N$, and $m_2$, the corresponding examples are sorted by the numbers
$N-20+1,\dots,N$.

Before making detailed comparisons of the obtained results
(s.~Tables 9--12
in
\ifthenelse{\boolean{AppendixOutsourcing}}
{\citet{HannesPaperAOutsourcing})}
{Appendix \ref{section:ResultTables})}
in
Subsubsections \ref{subsection:EasyExamplesWithNdividedBy2ConstraintComponents}--\ref{subsection:DifficultExamplesWithNConstraintComponents},
we give a short overview of them as a reason of clarity of the presentation: For all $N\in\lbrace20,40,60,80,100\rbrace$ we summarize the easy examples and the difficult examples, where we use two different line types for a better distinction of the comparisons of $m_2$
(for $m_2=\tfrac{N}{2}$ we use a dashed line
and for $m_2=N$ we use a solid line)
in Figures \ref{Figure:OverviewTableScalExEasy} and \ref{Figure:OverviewTableScalExDifficult},
which yields the following two summary tables consisting of the number of examples for which the reduced algorithm is better than MPBNGC resp.~SolvOpt (and vice versa) with respect to the solving time \TotalTime
\begin{center}
\begin{scriptsize}
\begin{tabular}{l|c|rrrrrrrr|}
\multicolumn{2}{c|}{(Color code: Grey)}
&
\multicolumn{4}{c|}{MPBNGC} & & \multicolumn{3}{|c|}{Red Alg}\\
\multicolumn{1}{c|}{}       &
\multicolumn{1}{c|}{}       &
\multicolumn{1}{c}{no termi-}     & 
\multicolumn{1}{c}{significantly} & \multicolumn{1}{c}{better} & \multicolumn{1}{c}{a bit} &
\multicolumn{1}{|c|}{nearly} &
\multicolumn{1}{c}{a bit} & \multicolumn{1}{c}{better} & \multicolumn{1}{c|}{significantly}
\\
\multicolumn{1}{c|}{Level}       &
\multicolumn{1}{c|}{$m_2$}       &
\multicolumn{1}{c}{nation}     & 
\multicolumn{1}{c}{better} & \multicolumn{1}{c}{} & \multicolumn{1}{c}{better} &
\multicolumn{1}{|c|}{equal} &
\multicolumn{1}{c}{better} & \multicolumn{1}{c}{} & \multicolumn{1}{c|}{better}
\\
\hline
Easy      & $\tfrac{N}{2}$ &   1 & 18 &	17 & 18 &	27 & 6 & 5 &	8\\
          & $N$            &   2 &  8 &	26 & 26 &	20 & 8 & 5 &  5\\
\hline
Difficult & $\tfrac{N}{2}$ &  73 &	0 &	 4 &  5 &  4 & 3 & 2 &  9\\
          & $N$            &  78 &	0 &	 1 &  1 &  5 & 1 & 4 & 10\\
\hline
\end{tabular}
\end{scriptsize}
\end{center}
\begin{center}
\begin{scriptsize}
\begin{tabular}{l|c|rrrrrrrr|}
\multicolumn{2}{c|}{(Color code: Black)}
&
\multicolumn{4}{c|}{SolvOpt} & & \multicolumn{3}{|c|}{Red Alg}\\
\multicolumn{1}{c|}{}       &
\multicolumn{1}{c|}{}       &
\multicolumn{1}{c}{no termi-}     & 
\multicolumn{1}{c}{significantly} & \multicolumn{1}{c}{better} & \multicolumn{1}{c}{a bit} &
\multicolumn{1}{|c|}{nearly} &
\multicolumn{1}{c}{a bit} & \multicolumn{1}{c}{better} & \multicolumn{1}{c|}{significantly}
\\
\multicolumn{1}{c|}{Level}       &
\multicolumn{1}{c|}{$m_2$}       &
\multicolumn{1}{c}{nation}     & 
\multicolumn{1}{c}{better} & \multicolumn{1}{c}{} & \multicolumn{1}{c}{better} &
\multicolumn{1}{|c|}{equal} &
\multicolumn{1}{c}{better} & \multicolumn{1}{c}{} & \multicolumn{1}{c|}{better}
\\
\hline
Easy      & $\tfrac{N}{2}$ & 18 & 14 &	25 & 11 &	15 &  6 &  7 &  4\\
          & $N$            & 11 & 16 &	21 & 11 &	15 & 10 &	11 &  5\\
\hline
Difficult & $\tfrac{N}{2}$ &  3 &  4 & 16 &  3 &	 8 & 15 &	28 & 23\\
          & $N$            &  0 &  5 &	 8 & 11 &	15 &  7 &	34 & 20\\
\hline
\end{tabular}
\end{scriptsize}
\end{center}
that are visualized in Figure \ref{Figure:OverviewTableScalExEasy} and Figure \ref{Figure:OverviewTableScalExDifficult}
\\
\begin{center}
\captionsetup{type=figure}
\includegraphics[width=8cm]{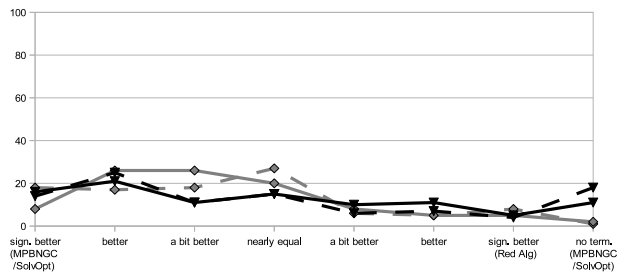}
\captionof{figure}{Easy examples (summary)}
\label{Figure:OverviewTableScalExEasy}
\end{center}
\begin{center}
\captionsetup{type=figure}
\includegraphics[width=8cm]{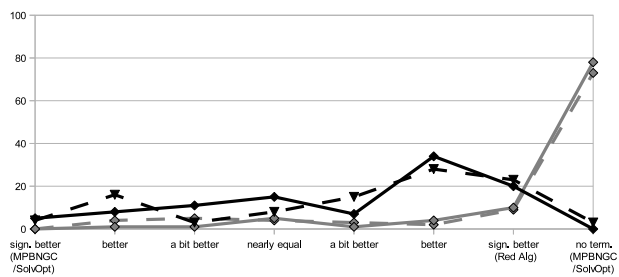}
\captionof{figure}{Difficult examples (summary)}
\label{Figure:OverviewTableScalExDifficult}
\end{center}
and that let us together with Figure \ref{NumericalResults:Figure:HigherDimensionalExamplesAllInOne}, in which the 
solving times $t_1$ for all examples are plotted
\begin{center}
\captionsetup{type=figure}
\includegraphics[width=16cm]{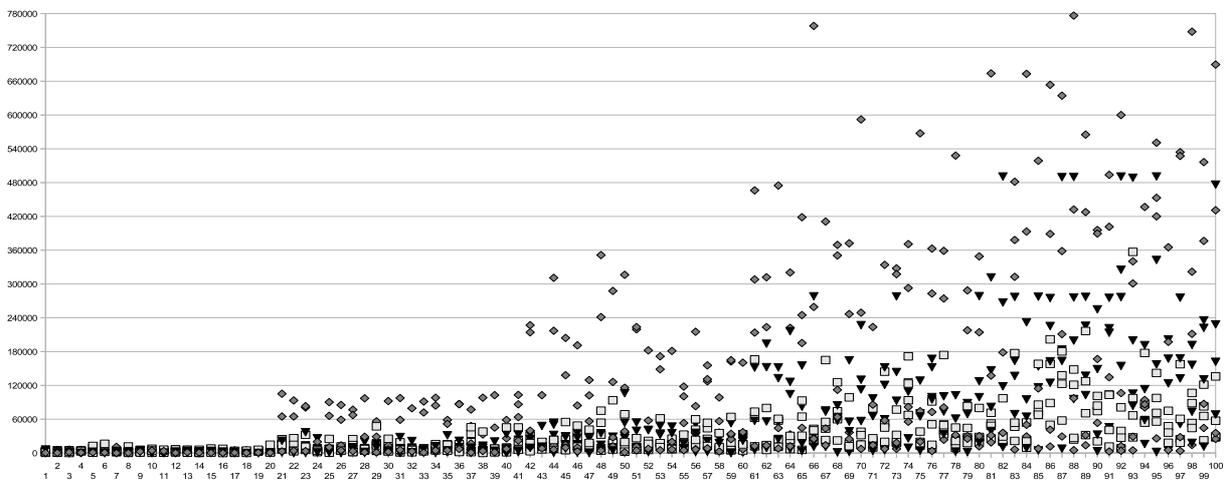}
\captionof{figure}{Solving time \TotalTime for all higher dimensional piecewise quadratic examples}
\label{NumericalResults:Figure:HigherDimensionalExamplesAllInOne}
\end{center}
draw the following conclusions:

For the easy examples the reduced algorithm is superior to MPBNGC in thirteen percent of the examples, for thirty percent of the examples one of these two solvers has small advantages over the other (in most cases MPBNGC is the slightly more successful one), the performance differences between the considered two algorithms can be completely neglected for one quarter of the examples, and for one third of the examples MPBNGC beats the reduced algorithm clearly.
MPBNGC was not able to terminate successfully for many of the difficult examples in particular for $N\in\lbrace60,80,100\rbrace$ despite significantly longer running times as it can be seen in Figure \ref{NumericalResults:Figure:HigherDimensionalExamplesAllInOne} (in additional test runs with a softer termination criterion MPBNGC did terminate for approximately half of the difficult examples, but the quality of the obtained minimizers was not comparable with the corresponding $\hat{x}$ produced by the reduced algorithm, while for the comparisons presented here this quality is comparable) and therefore the reduced algorithm is superior to MPBNGC in 88 percent of these examples. Furthermore, for five percent of the examples one of these two solvers has small advantages over the other, the performance differences between the considered two algorithms can be completely neglected for further five percent of the examples, and for the remaining two percent of the examples MPBNGC beats the reduced algorithm clearly.

For the easy examples the reduced algorithm is superior to SolvOpt in thirty percent of the examples, for fifteen percent of the examples one of these two solvers has small advantages over the other, the performance differences between the considered two algorithms can be completely neglected for further fifteen percent of the examples, and in the remaining forty percent of the examples SolvOpt beats the reduced algorithm clearly.
For the difficult examples the reduced algorithm is superior to SolvOpt in a bit more than half of the examples (including many examples with $N\in\lbrace80,100\rbrace$), for twenty percent of the examples one of these two solvers has small advantages over the other, the performance differences between the considered two algorithms can be completely neglected for ten percent of the examples, and in the remaining (a bit less than) twenty percent of the examples SolvOpt beats the reduced algorithm clearly.
In particular note that only very few $F_j$ are active at the points which SolvOpt found at termination for the easy examples (in comparison to both the reduced algorithm and MPBNGC), which might indicate that SolvOpt has some problems coming very close to the boundary. Although this behavior improves for the difficult examples, there still remains a clear gap in the number of active $F_j$ between SolvOpt and the other two solvers.


We want to emphasize the reduced algorithm was the only solver that terminated for all higher dimensional examples successfully, i.e.~with a stationary point that is sufficiently accurate. Moreover, the solving times of the reduced algorithm are quite stable over all dimensions $N\in\lbrace20,40,60,80,100\rbrace$.
\begin{remark}
Since MOSEK supports multiple CPUs in particular for solving
QCQPs (cf.~\citet[p.152,~8.1.4~Using~multiple~CPU's]{AndersenMosek}), we may expect faster solving times for the reduced algorithm on such a system in particular for higher dimensional problems. Nevertheless, we have not been able to test this yet.

We also expect a significant improvement of the full algorithm if a QCQP-solver is used which exploits the special structure of the QCQP (\ref{BundleSQPmitQCQP:Alg:QPTeilproblem}).
\end{remark}

\subsubsection{Easy examples with N/2 constraint components}
\label{subsection:EasyExamplesWithNdividedBy2ConstraintComponents}
We summarize the investigations of the results of the easy examples with $m_2:=\tfrac{N}{2}$, which can be found in Table 9
in
\ifthenelse{\boolean{AppendixOutsourcing}}
{\citet{HannesPaperAOutsourcing}}
{Appendix \ref{section:ResultTables}}
and which are visualized in Figure 30 in \citet{HannesPaperAOutsourcing}.

\paragraph{Reduced algorithm vs. MPBNGC}
MPBNGC does not satisfy its termination criterion for one example
within the given number of iterations and function evaluations. For the remaining 99 examples we obtain that the reduced algorithm is
significantly better in 8 examples,
better in 5 examples
and
a bit better in 6 examples
in comparison with MPBNGC, while MPBNGC is
significantly better in 18 examples,
better in 17 examples
a bit better in 18 examples
in comparison with the reduced algorithm, and in 27 examples the solving times of both algorithms do not differ significantly.

\paragraph{Reduced algorithm vs. SolvOpt}
SolvOpt does not satisfy its termination criterion for 18 examples
within the given number of iterations and function evaluations. For the remaining 82 examples we obtain that the reduced algorithm is
significantly better in 4 examples,
better in 7 examples
and
a bit better in 6 examples
in comparison with SolvOpt, while SolvOpt is
significantly better in 14 examples,
better in 25 examples
and
a bit better in 11 examples
in comparison with the reduced algorithm, and in 15 examples the solving times of both algorithms do not differ significantly.

\subsubsection{Easy examples with N constraint components}
\label{subsection:EasyExamplesWithNConstraintComponents}
We summarize the investigations of the results of the easy examples with $m_2:=N$, which can be found in Table 10
in
\ifthenelse{\boolean{AppendixOutsourcing}}
{\citet{HannesPaperAOutsourcing}}
{Appendix \ref{section:ResultTables}}
and which are visualized in Figure 31 in \citet{HannesPaperAOutsourcing}.

\paragraph{Reduced algorithm vs. MPBNGC}
MPBNGC does not satisfy its termination criterion for two examples
within the given number of iterations and function evaluations. For the remaining 98 examples we obtain that the reduced algorithm is
significantly better in 5 examples,
better in 5 examples
and
a bit better in 8 examples
in comparison with MPBNGC, while MPBNGC is
significantly better in 8 examples,
better in 26 examples
and
a bit better in 26 examples
in comparison with the reduced algorithm, and in 20 examples the solving times of both algorithms do not differ significantly.

\paragraph{Reduced algorithm vs. SolvOpt}
SolvOpt does not satisfy its termination criterion for 11 examples
within the given number of iterations and function evaluations. For the remaining 89 examples we obtain that the reduced algorithm is
significantly better in 5 examples,
better in 11 examples
and
a bit better in 10 examples
in comparison with SolvOpt, while SolvOpt is
significantly better in 16 examples,
better in 21 examples
and
a bit better in 11 examples
in comparison with the reduced algorithm, and in 15 examples the solving times of both algorithms do not differ significantly.

\subsubsection{Difficult examples with N/2 constraint components}
\label{subsection:DifficultExamplesWithNdividedby2ConstraintComponents}
We summarize the investigations of the results of the difficult examples with $m_2:=\tfrac{N}{2}$, which can be found in Table 11
in
\ifthenelse{\boolean{AppendixOutsourcing}}
{\citet{HannesPaperAOutsourcing}}
{Appendix \ref{section:ResultTables}}
and which are visualized in Figure 32 in \citet{HannesPaperAOutsourcing}.

\paragraph{Reduced algorithm vs. MPBNGC}
MPBNGC does not satisfy its termination criterion for 73 examples
within the given number of iterations and function evaluations. For the remaining 27 examples we obtain that the reduced algorithm is
significantly better in 9 examples,
better in 2 examples
and
a bit better in 3 examples
in comparison with MPBNGC, while MPBNGC is
better in 4 examples
and
a bit better in 5 examples
in comparison with the reduced algorithm, and in 4 examples the solving times of both algorithms do not differ significantly.

\paragraph{Reduced algorithm vs. SolvOpt}
SolvOpt does not satisfy its termination criterion for 3 examples
within the given number of iterations and function evaluations. For the remaining 97 examples we obtain that the reduced algorithm is
significantly better in 23 examples,
better in 28 examples
and
a bit better in 15 examples
in comparison with SolvOpt, while SolvOpt is
significantly better in 4 examples,
better in 16 examples
and
a bit better in 3 examples
in comparison with the reduced algorithm, and in 8 examples the solving times of both algorithms do not differ significantly.

\subsubsection{Difficult examples with N constraint components}
\label{subsection:DifficultExamplesWithNConstraintComponents}
We summarize the investigations of the results of the difficult examples with $m_2:=N$, which can be found in Table 12
in
\ifthenelse{\boolean{AppendixOutsourcing}}
{\citet{HannesPaperAOutsourcing}}
{Appendix \ref{section:ResultTables}}
and which are visualized in Figure 33 in \citet{HannesPaperAOutsourcing}.

\paragraph{Reduced algorithm vs. MPBNGC}
MPBNGC does not satisfy its termination criterion for 78 examples
within the given number of iterations and function evaluations. For the remaining 22 examples we obtain that the reduced algorithm is
significantly better in 10 examples,
better in 4 examples
and
a bit better in 1 example
in comparison with MPBNGC, while MPBNGC is
better in 1 example
and
a bit better in 1 example
in comparison with the reduced algorithm, and in 5 examples the solving times of both algorithms do not differ significantly.

\paragraph{Reduced algorithm vs. SolvOpt}
For our 100 examples we obtain that the reduced algorithm is
significantly better in 20 examples,
better in 34 examples
and
a bit better in 7 examples
in comparison with SolvOpt, while SolvOpt is
significantly better in 5 examples,
better in 8 examples
and
a bit better in 11 examples
in comparison with the reduced algorithm, and in 15 examples the solving times of both algorithms do not differ significantly.

\section{Conclusion}
In this paper we investigated numerical aspects of the feasible second
order bundle algorithm for nonsmooth, nonconvex optimization problems
with inequality constraints. Since one of the main characteristics of
this method is that the search direction is determined by solving a
convex QCQP, we investigated certain versions of the search direction
problem and we justified the version chosen by us numerically by
comparing the results of different solvers for the computation of the
search direction. Furthermore, we made comparisons between the test
results of our implementation of the second order bundle algorithm,
MPBNGC by \citet{MPBNGC} and SolvOpt by \citet{SolvOptPublication} for
some examples of the Hock-Schittkowski collection by
\citet{Schittkowski,Schittkowski2} and for custom examples that arise
in the context of finding exclusion boxes for quadratic CSPs, where
for both of these types of examples we were able to achieve good
results with respect to the number of evaluations of function
dependent data, as well as for higher dimensional piecewise quadratic
examples, in which our implementation achieved good results in
comparison with the other solvers in particular in the case that many
constraint components were active at the solution. Summarizing the
results it can be seen that the the SQP-like algorithm tends to
compare the better the higher the dimension of the problem and the
more difficult the nonsmoothness around the optimal point are.





\end{document}